\documentclass{elsarticle}[preprint,review,10pt]

\journal{Journal of Computational Physics}

\usepackage{graphicx} % Required for inserting images
\usepackage{algorithm}
\usepackage{algorithmic}

\usepackage[utf8]{inputenc}
\usepackage{graphicx}
\usepackage{bm}
\usepackage{comment}
\usepackage{psfrag}
\usepackage{latexsym,amsmath,amsfonts,amscd,amsthm}
\usepackage{changebar}
\usepackage{color}
\usepackage{bm}
\usepackage{tikz}
\usepackage{multirow}
\usepackage{xcolor}
\usepackage{hyperref}
\usepackage{setspace}
\usepackage{amssymb}
\usepackage{mathrsfs}
\usepackage{caption}
\usepackage{subcaption}
\usepackage[title]{appendix}
\usepackage{ulem}
\usepackage{tikz}
\usetikzlibrary{arrows,backgrounds}

\newlength{\spse}
\newtheorem{thm}{Theorem}[section]

\newtheorem{rem}[thm]{Remark}

\newcommand{\BU}{\mathbf{U}}

\newcommand{\BT}{\mathbf{T}}
\newcommand{\BI}{\mathbf{I}}

\newcommand{\BA}{\mathbf{A}}
\newcommand{\WBA}{\widetilde{\mathbf{A}}}
\newcommand{\BQ}{\mathbf{Q}}

\newcommand{\BOmega}{\boldsymbol{\Omega}}

\newcommand{\bx}{\mathbf{x}}
\newcommand{\bu}{\mathbf{u}}
\newcommand{\bw}{\mathbf{w}}
\newcommand{\bc}{\mathbf{c}}
\newcommand{\bb}{\mathbf{b}}
\newcommand{\wbb}{\widetilde{\mathbf{b}}}
\newcommand{\bv}{\boldsymbol{v}}

\newcommand{\bn}{\boldsymbol{n}}
\newcommand{\bxi}{\boldsymbol{\xi}}

\newcommand{\half}{\frac{1}{2}}
\newcommand{\brho}{\boldsymbol{\rho}}
\newcommand{\bg}{\boldsymbol{g}}
\newcommand{\br}{\boldsymbol{r}}

\newcommand{\bG}{\boldsymbol{G}}
\newcommand{\wbG}{\widetilde{\boldsymbol{G}}}

\newcommand{\mfw}{\mathfrak{w}}

\newcommand{\bpsi}{{{\boldsymbol{\psi}}}}
\newcommand{\bphi}{{{\boldsymbol{\phi}}}}

\newcommand{\BJ}{{\bf{J}} }
\newcommand{\BC}{{\bf{C}} }
\newcommand{\BH}{{\bf{H}} }
\newcommand{\BM}{{\bf{M}} }
\newcommand{\BD}{{\bf{D}}}
\newcommand{\BZ}{{\bf{Z}}}
\newcommand{\BSigma}{{\mathbf{\Sigma}}}
\newcommand{\bmu}{{\boldsymbol{\mu}}}
\newcommand{\be}{{\bf{e}}}
\newcommand{\bz}{{\bf{z}}}
\newcommand{\by}{{\bf{y}}}
\newcommand{\boldeta}{{\boldsymbol{\eta}}}

\newcommand{\pzc}{\textcolor{black}}

\newcommand{\epsrom}{\epsilon_{\textrm{ROM}}}

\newcommand{\BS}{\mathbf{S}}

\title{A Flexible GMRES Solver with Reduced Order Model Enhanced  Synthetic Acceleration Preconditioenr for Parametric Radiative Transfer Equation}

\begin{document}

\author{
Zhichao Peng
}

\ead{pengzhic@ust.hk}

\affiliation{organization={Department of Mathematics, The Hong Kong University of Science and Technology},%Department and Organization
            addressline={Clear Water
Bay, Kowloon},
city={Hong Kong},
country={China}
}

%%%%%%%%%%%%%%%%%%%%%%%%%%%%%%%%%%%%%%%%%%%%%%%%%%%%%%%%%%%%
\begin{abstract}
Parametric  radiative transfer equation (RTE) occurs in multi-query applications such as uncertainty quantification, inverse problems, and sensitivity analysis, which require solving RTE multiple times for a range of parameters. Consequently, efficient iterative solvers are highly desired.

Classical Synthetic Acceleration (SA) preconditioners for RTE build on  low order approximations to an ideal kinetic correction equation such as its diffusion limit in Diffusion Synthetic Acceleration (DSA). Their performance depends on the effectiveness of the underlying low order approximation. In addition, they do not leverage low rank structures with respect to the parameters of the parametric problem.

To address these issues, we proposed a ROM-enhanced SA strategy, called ROMSAD, under the Source Iteration framework in Peng (2024).  In this paper, we further extend the ROMSAD preconditioner to flexible general minimal residual method (FGMRES). The main new advancement is twofold. First, after identifying the ideal kinetic correction equation within the FGMRES framework, we reformulate it into an equivalent form, allowing us to develop an iterative procedure to construct a ROM for this ideal correction equation without directly solving it. Second, we introduce a greedy algorithm to build the underlying ROM for the ROMSAD preconditioner more efficiently.

Our numerical examples demonstrate that FGMRES with the ROMSAD preconditioner (FGMRES-ROMSAD) is more efficient than GMRES with the right DSA preconditioner. Furthermore, when the underlying ROM in ROMSAD is not highly accurate, FGMRES-ROMSAD exhibits greater robustness compared to Source Iteration accelerated by ROMSAD.
\end{abstract}


\maketitle

%%%%%%%%%%%%%%%%%%%%%%%%%%%%%%%%%%%%%%%%%%%%%%%%%%%%%%%%%%%%
\section{Introduction}
Radiative transfer equation (RTE) provides models for particle systems in a wide range of applications, e.g., medical imaging \cite{arridge2009optical}, remote sensing \cite{spurr2001linearized}, nuclear engineering \cite{pomraning1973equations}, and astrophysics \cite{janka2007theory}. In multi-query engineering or scientific applications, such as design optimization, inverse problems, sensitivity analysis, and uncertainty quantification, RTE is often solved repeatedly within a parametric setting.  The parameters for these parametric problems typically originate from the parametrization of material properties, boundary conditions, or geometric configurations. Due to the need for solving the RTE multiple times, efficient numerical solvers for parametric RTE are highly desired.

Among iterative linear solvers for RTE, Source Iteration (SI) and Krylov solvers with Synthetic Acceleration (SA) \cite{adams2001discontinuous,warsa2004krylov,azmy2010advances} are widely used. It is well known that \pzc{classical} SI may converge arbitrarily slowly \cite{adams2001discontinuous}. To accelerate its convergence, SA preconditioner introduces a correction to the scalar flux (also known as the macroscopic density). The ideal scalar flux correction, which guarantees the convergence of SI in the next iteration, is the integral of an ideal angular flux correction, which is the solution to an ideal kinetic correction equation. However, solving this ideal kinetic correction equation results in almost the same computational costs as directly solving the original problem. Hence, in practice, a computationally cheap low-order approximation to it is employed. Popular choices include the diffusion limit of the kinetic correction equation in Diffusion Synthetic Acceleration (DSA) \cite{kopp1963synthetic,alcouffe1977dittusion,adams1992diffusion,wareing1993new,Adams2002FastIM}, the variable Eddington factor in
Quasi-Diffusion method \cite{gol1964quasi,anistratov1993nonlinear,olivier2023family}, and the $S_2$ approximation in S2SA \cite{lorence1989s}. However, as pointed out in \cite{azmy2002unconditionally,warsa2004krylov}, the effectiveness of SI-SA  may degenerate for challenging problems, such as those involving material discontinuities. To overcome this degeneration, \cite{warsa2004krylov} extends SA preconditioners to Krylov solvers by rewriting the SI-SA as a left preconditioend Richardson iteration for a memory-efficient discrete formulation of the RTE.

Despite the great success and popularity of these classical SA strategies, their effectiveness hinges on the underlying empirical low order approximations. For instance, when the problem is far from its diffusion limit, DSA may lose its efficiency \cite{ren2019fast}.  Moreover, classical SA strategies also do not leverage low rank structures with respect to parameters of parametric problems. 

To tackle these issues, data-driven reduced-order models (ROMs) \cite{benner2015survey} can be utilized. ROM is a technique to build low rank approximations leveraging low rank structures in parametric problems. Recently, data-driven ROMs have been developed for RTE in \cite{buchan2015pod, tencer2017accelerated, choi2021space, tano2021affine, behne2022minimally, behne2023parametric, halvic2023non, mcclarren2022data, peng2022reduced, coale2023reduced, coale2024reduced, peng2024micro, buchan2024reduced, hardy2024proper,matsuda2025rbmrte}. In \cite{peng2024romsad}, under the SI framework, we employed a data-driven ROM for the ideal kinetic correction equation to design a ROM-enhanced SA strategy called ROMSAD  to leverage the original kinetic description of the correction equation and low rank structures in parametric problems. The ideal kinetic correction equation  varies as source iterations continue.  To take this iteration dependence into account without consuming excessive memory, we built the ROM using the ideal corrections for the first few iterations. Without data for later iterations, the accuracy of the ROM may reduce as source iterations continue. Hence, in our ROMSAD method, we applied ROM-based corrections in the first few iterations to improve its efficiency  and then switched to DSA to maintain robustness in later iterations. Using sufficiently accurate ROM, SI-ROMSAD is able to achieve significant acceleration over SI-DSA. 

However, the ROMSAD preconditioner is only built for SI in \cite{peng2024romsad}, while Krylov solvers may be preferred in many applications for their faster convergence and better robustness. Additionally, the underlying ROM is constructed using high-fidelity ideal corrections for parameters across the entire training set, leading to potentially unnecessary computational costs due to the generation of all this high-fidelity data. In this paper, we extend our ROMSAD preconditioner from SI to the Krylov framework and apply a greedy algorithm \cite{hesthaven2016certified,matsuda2025rbmrte} to build the underlying ROM  more efficiently. 
\begin{enumerate}
    \item To extend our method to the Krylov solver, following the extension of classical SA in \cite{warsa2004krylov}, we first rewrite ROMSAD preconditioned SI as a Richardson iteration with a nonlinear left preconditioner solving a memory-efficient discrete formulation. Next, we identify the ideal kinetic correction equation in Krylov methods and reformulate it to efficiently obtain its solution without directly solving it. Due to the switching step in ROMSAD, the resulting preconditioner is nonlinear. Therefore, we choose the flexible general minimal residual method (FGMRES) \cite{saad1993flexible} as the underlying Krylov solver.
    \item When constructing the ROM, we use a greedy algorithm to iteratively enrich the underlying reduced order space by adding corrections for the most ``representative" parameters. These ``representative" parameters are identified by a residual-based error indicator. This greedy algorithm improves the offline efficiency by  generating high-fidelity data only for the identified ``representative" parameters rather than the entire training set.
\end{enumerate}
Compared to our previous work \cite{peng2024romsad},
the main advancements are as follows. First, we   reformulate the ideal kinetic correction equation for FGMRES to allow efficiently obtaining its solution without solving it. Second, we apply a greedy algorithm to improve the efficiency of the ROM construction. 

Our numerical tests demonstrate that FGMRES with the ROMSAD preconditioner (FGMRES-ROMSAD) is more efficient than right preconditioned GMRES using the classical DSA preconditioner. Furthermore, when the underlying ROM is not highly accurate, FGMRES-ROMSAD  is more robust than SI accelerated by ROMSAD.

Before detailing our method, we first contextualize it by briefly reviewing data-driven ROM-based accelerations of iteartive linear solvers for RTE. Dynamic Mode Decomposition (DMD) is exploited as a low rank update strategy in SI \cite{mcclarren2022data} and eigenvalue solvers in \cite{mcclarren2019calculating,roberts2019acceleration,smith2023variable}.  A neural network surrogate for the transport sweep in SI is developed in \cite{tano2021sweep}; however, this surrogate is only able to provide low-fidelity approximations.  Random Singular Value Decomposition has been applied to build a low rank boundary-to-boundary map to enhance a Schwartz solver for RTE \cite{chen2021low}. A ROM-based adaptive nested iterative solver is proposed for the Discontinuous Petrov Galerkin finite element method solving RTE \cite{dahmen2020adaptive}. Besides accelerations based on data-driven low rank approximations for RTE, Tailored Finite Point Method (TFPM) are utilized in   \cite{fu2024fast,song2024adaptive}, while low rank matrix or tensor decomposition techniques are leveraged in the time marching \cite{peng2020low,einkemmer2021asymptotic,peng2021high,kusch2021robust,yin2024towards,einkemmer2024asymptotic,sands2024high} and iterative solvers of RTE \cite{bachmayr2024low}.

 Beyond the scope of RTE, ROM-based preconditioners have been developed for elliptic and parabolic problems in \cite{pasetto2017reduced, santo2018multi, cortes2018pod, diaz2021accelerating, hou2023reduced}. The ROM-based two-level preconditioner for elliptic equations in \cite{santo2018multi} shares similarities with the ROMSAD method, as both methods build single or multiple ROMs using data from multiple iterations. We will show that the projection step to construct the ROM-based preconditioner in these two methods has a major difference due to their distinct starting points. Our method is motivated by SA based on the ideal kinetic correction equation, while the method in \cite{santo2018multi} is inspired by two-grid preconditioners.

This paper is organized as follows. We briefly review the model equation and discretizations for it in Sec. \ref{sec:background}, followed by basic ideas of classical SA in Sec. \ref{sec:sa-review}. Then, we outline how to extend ROMSAD method from SI to Krylov method and the greedy algorithm to build ROMs more efficiently  in Sec. \ref{sec:romsad-preconditioner}. The proposed method will be numerically tested in Sec. \ref{sec:numerical}. At last, we draw our conclusions in Sec. \ref{sec:conclusion}.
%%%%%%%%%%%%%%%%%%%%%%%%%%%%%%%%%%%%%%%%%%%%%%%%%%%%%%%%%%%%
\section{Model equation and full order discretization\label{sec:background}}
In this paper, we consider steady state linear RTE with one energy group, isotropic scattering, source and isotropic inflow boundary conditions \pzc{defined on the computational domain $\Gamma_{\bx}$ in the physical space and the angular space $\mathbb{S}^2$}:
\begin{subequations}
\label{eq:rte}
    \begin{align}
    &\BOmega \cdot \pzc{\nabla} \psi(\bx,\BOmega) + \sigma_t(\bx) \psi(\bx,\BOmega) = \sigma_s(\bx) \phi(\bx) + G(\bx), \quad\sigma_t(\bx)=\sigma_s(\bx)+\sigma_a(\bx),\\
    &\phi(\bx) = \frac{1}{4\pi}\int_{\mathbb{S}^2} \psi(\bx,\BOmega) d\BOmega, \quad \bx\in\Gamma_{\bx},\\
    &\psi(\bx,\BOmega) = g(\bx), \quad \bx\in \partial\Gamma_{\bx},\;\BOmega\cdot \bn(\bx)<0.\label{eq:boundary-condition}
    \end{align}
\end{subequations}
Here, $\psi(\bx,\BOmega)$ is the angular flux (also known as radiation intensity or particle distribution) corresponding to angular direction on the unit sphere $\BOmega=(\BOmega_x,\BOmega_y,\BOmega_z)\in\mathbb{S}^2$ and spatial location $\bx\in\Gamma_{\bx}$,  $\phi(\bx)=\frac{1}{4\pi}\int_{\mathbb{S}^2} \psi(\bx,\BOmega) d\BOmega$ is the scalar flux (also known as the macroscopic density), $G(\bx)$ is an isotropic source, $\sigma_t(\bx)=\sigma_s(\bx)+\sigma_a(\bx)\geq 0$ is the total cross section, $\sigma_s(\bx)\geq 0$ is  the isotropic scattering cross section, and $\sigma_a(\bx)\geq 0$ is the absorption cross section. In the inflow boundary condition \eqref{eq:boundary-condition}, $\bn(\bx)$ is the outward normal direction of the computational domain $\Gamma_{\bx}$ at location $\bx\in\partial\Gamma_{\bx}$.
\pzc{For simplicity, in this paper, we focus on the case where the spatial dimension is two and the angular space is the unit sphere $\mathbb{S}^2$.}

As the scattering effect goes to infinity, i.e. $\sigma_s(\bx)\rightarrow\infty$, RTE \eqref{eq:rte} converges to its diffusion limit: $\psi(\bx,\BOmega)\rightarrow\phi(\bx)$ satisfying
\begin{equation}
-\nabla_x\cdot\left(\mathcal{D}\frac{1}{\sigma_s}\nabla\phi\right)=-\sigma_a\phi+G,\quad 
\text{where}\quad\mathcal{D}=\frac{1}{4\pi}\textrm{diag}(\int_{\mathbb{S}^2}\BOmega_x^2d\BOmega,\int_{\mathbb{S}^2}\BOmega_y^2d\BOmega,\int_{\mathbb{S}^2}\BOmega_z^2d\BOmega).
    \label{eq:diffusion-limit}
\end{equation}
When discretizing RTE  \eqref{eq:rte}, this diffusion limit \pzc{needs} to be preserved on the discrete level. In other words, an asymptotic preserving method \cite{jin2010asymptotic} should be applied  
to capture this diffusion limit without using a small mesh size $h$ to resolve the mean free path of particles: $\sigma_t h\ll 1$. In this paper, we employ the discrete ordinates ($S_N$) angular discretization and upwind discontinuous Galerkin (DG) spatial discretization, which is proved to be asymptotic preserving \cite{adams2001discontinuous,guermond2010asymptotic}.

%%%%%%%%%%%%%%%%%%%%%%%%%%%%%%%%%%%%%%%%%%%%%%%%%%%%%%%%%%%%
\subsection{Discrete ordinates ($S_N$) angular discretization}
The basic idea of the discrete ordinates ($S_N$) method \cite{pomraning1973equations} is to sample the angular flux at a set of quadrature points $\{\BOmega_j\}_{j=1}^{N_{\BOmega}}$ with quadrature weights $\{\pzc{\omega_j}\}_{j=1}^{N_{\BOmega}}$ in the angular space. The scalar flux is approximated by  the underlying quadrature rule. Here, we use a normalized  quadrature rule 
$\{(\BOmega_j,\omega_j)\}_{j=1}^{N_{\BOmega}}$ satisfying 
$$
\sum_{j=1}^{N_{\BOmega}}\omega_j=1\quad\text{and}\quad\frac{1}{4\pi}\int_{\mathbb{S}^2}g(\BOmega)d\BOmega\approx \sum_{j=1}^{N_{\BOmega}}\omega_j g(\BOmega_j).
$$
Then, the angular discretized version of \eqref{eq:rte} is given by:
\begin{subequations}
\label{eq:sn_system}
\begin{align}
&(\BOmega_j\cdot\nabla_x+\sigma_t(\bx) )\psi(\bx,\BOmega_j) = \sigma_s(\bx) \phi(\bx)+G(\bx),\quad \phi(\bx)= \sum_{j=1}^{N_{\BOmega}} \omega_j \psi(\bx,\BOmega_j),\\
&\psi(\bx,\BOmega_j) = g(\bx), \quad \bx\in\partial\Gamma_{\bx},\quad\BOmega_j\cdot\bn(\bx)<0.
\end{align}
\end{subequations}
In this paper, we use Chebyshev-Legendre (CL) quadrature rule. The CL quadrature rule CL($N_{\theta},N_{\BOmega_z}$) is  the tensor product of the normalized $N_{\theta}$-\pzc{points} Chebyshev quadrature rule for the unit circle 
\begin{equation}
\left\{(\theta_j,\omega_{\theta,j}): \;\theta_j = \frac{2j\pi}{N_\theta}-\frac{\pi}{N_\theta}\;\text{and}\;\omega^\phi_j=\frac{1}{N_\theta},\; j=1,\dots,N_\theta\right\}
\end{equation}
and the normalized $N_{\BOmega_z}$-points Gauss-Legendre quadrature rule $\{(\BOmega_{z,j},\omega_{z,j})\}_{j=1}^{N_{\BOmega_z}}$ for the $z$-component of angular direction  $\BOmega_z\in[-1,1]$ satisfying $\sum_{j=1}^{N_{\BOmega_z}} \omega_{z,j} =1$. 
The quadrature points and the associated quadrature weights of the CL($N_{\theta}$,$N_{\BOmega_z}$) quadrature rule are defined as
\begin{equation}
\BOmega_j=\left(\cos(\theta_{j_1})\sqrt{1-\BOmega_{z,j_2}^2},\sin(\theta_{j_1})\sqrt{1-\BOmega_{z,j_2}^2},\BOmega_{z,j_2}\right)\quad\text{and}\quad\omega_j=\omega_{\theta,j_1}\omega_{\BOmega_z,j_2},
\end{equation}
where  $1\leq j_1\leq N_\theta$, $1\leq j_2\leq N_{\BOmega_z}$ and $1\leq j=(j_2-1) N_{\theta}+j_1\leq \pzc{N_{\BOmega}}$, where \pzc{$N_{\BOmega}=N_{\BOmega_z}N_{\theta}$.}

%%%%%%%%%%%%%%%%%%%%%%%%%%%%%%%%%%%%%%%%%%%%%%%%%%%%%%%%%%%%
\subsection{Upwind discontinuous Galerkin spatial discretization}

For simplicity, we consider a  2D rectangular computational domain $\Gamma_{\bx}=[x_l,x_r]\times[y_l,y_r]$ and partition it with rectangular meshes $\mathcal{T}_h=\{\mathcal{T}_i,\;\forall \mathcal{T}_i\;\text{being rectangle}\}_{i=1}^{N_{\bx}}$. Denote the set of cell edges as $\partial \mathcal{T}_h$ and the set of edges on the inflow boundary for $\BOmega_j$ as 
\begin{equation*}
\partial \mathcal{T}_{h,j}^{(\textrm{ibc})}=\{\mathcal{E}: \mathcal{E}\in\partial\mathcal{T}_h,\;\mathcal{E}\subset\partial\Omega_{\bx},\;\BOmega_j\cdot\bn(\bx)<0,\forall \bx\in\mathcal{E}\},
\end{equation*}
where $\bn(\bx)$ is the outward normal direction of $\Gamma_{\bx}$ at $\bx$.

To discretize the discrete ordinate equation \eqref{eq:sn_system} in space, we use a $Q^K$ upwind discontinuous Galerkin (DG) method which seeks solution in the finite element space
\begin{equation}
       U_h^K(\mathcal{T}_h):=\{u(\bx): u(\bx)|_{\mathcal{T}_{i}}\in Q^K(\mathcal{T}_{i}),1\leq i \leq N_{\bx}\}.\label{eq:discrete_space}
\end{equation}
Here, $Q^K(\mathcal{T}_{i})$ is the space of bi-variate polynomials on the element $\mathcal{T}_{i}$ whose degree in each direction is
at most $K$. 
Specifically, we seek $\psi_h(\bx,\BOmega_j)\in U_h^K(\mathcal{T}_h)$, $j=1,\dots,N_{\BOmega}$ satisfying $\forall \eta_h(\bx)\in U_h^K(\mathcal{T}_h),$
\begin{align}
-\sum_{i=1}^{N_{\bx}}&\int_{\mathcal{T}_i} \Big(\BOmega_j\cdot\nabla\eta_h(\bx)\Big) \psi_h(\bx,\BOmega_j) d\bx+\sum_{\mathcal{E}\in\partial\mathcal{T}_h\setminus\partial\mathcal{T}_{h,j}^{(\textrm{ibc})}} \int_{\mathcal{E}} \widehat{\BH}(\BOmega_j,\psi_h, \bn(\bx))\eta_h(\bx) d\bx+\sum_{i=1}^{N_{\bx}}\int_{\mathcal{T}_i}\sigma_t(\bx) \psi_h(\bx,\BOmega_j)\eta_h(\bx) d\bx
\notag
\\
= &\sum_{i=1}^{N_{\bx}}\int_{\mathcal{T}_i}\sigma_s(\bx) \phi_h(\bx)\eta_h(\bx) d\bx
 + \sum_{i=1}^{N_{\bx}}\int_{\mathcal{T}_i}G(\bx)\eta_h(\bx) d\bx
 -\sum_{\mathcal{E}\in\partial\mathcal{T}_{h,j}^{(\textrm{ibc})}}\int_{\mathcal{E}} g(\bx,\BOmega_j)\phi_h(\bx) \BOmega_j\cdot \bn(\bx) d\bx.
 \label{eq:DG}
\end{align}
Here, $\phi_h(\bx)=\sum_{j=1}^{N_{\BOmega}}\omega_j\psi_h(\bx,\BOmega_j)$.
The upwind numerical flux $\widehat{\BH}(\BOmega_j, \psi_h,\bn(\bx))$ along the edge $\mathcal{E}$ shared by an element $\mathcal{T}^-$ and its neighbor $\mathcal{T}^+$ is defined as
\begin{align}
\widehat{\BH}(\BOmega_j, \psi_h,\bn(\bx))\Big|_{\mathcal{E}} = \frac{\BOmega_j\cdot\bn(\bx)}{2}\Big(\psi_h^+(\bx,\BOmega_j)+\psi_h^-(\bx,\BOmega_j)\Big)+\frac{|\BOmega_j\cdot \bn(\bx)|}{2}\Big(\psi_h^-(\bx,\BOmega_j)-\psi_h^+(\bx,\BOmega_j)\Big),
\label{eq:upwind}
\end{align}
where $\psi_h^{\pm}(\bx,\BOmega_j)$ is the restriction of $\psi_h(\bx,\BOmega_j)$ onto $\mathcal{T}^{\pm}$, and $\bn(\bx)$ is  the unit outward normal direction with respect to the element $\mathcal{T}^-$. 
When the polynomial order $K\geq 1$, the above discrete ordinate DG method is asymptotic preserving \cite{adams2001discontinuous,guermond2010asymptotic}. 

\subsection{Matrix-vector formulation}
Assume $\{\eta_i(\bx)\}_{i=1}^{N_{\textrm{DOF}}}$ be an orthonormal basis for $U_h^K(\mathcal{T}_h)$. The degrees of freedom $\bpsi_j$ for $\psi_h(\bx,\BOmega_j)$ satisfies 
\begin{equation}
\psi_h(\bx,\BOmega_j)=\bpsi_j^T(\eta_1(\bx),\dots,\eta_{N_{\textrm{DOF}}}(\bx))^T =\sum_{k=1}^{N_{\textrm{DOF}}}\bpsi_{j,k}\eta_k(\bx).
\end{equation}
The matrix-vector form of the fully discrete scheme is: 
\begin{subequations}
\label{eq:dg_matrix_vec}
\begin{align}
    (\BD_j+\BSigma_t) \bpsi_j = \BSigma_s \BQ\bpsi+\bG+\bg_j^{(\textrm{bc})}=\BSigma_s\bphi+\bG+\bg_j^{(\textrm{bc})}=\BSigma_s\bphi+\widetilde{\bG}_j,\label{eq:discrete-rte}\\
    \bpsi=(\bpsi_1,\dots,\bpsi_{N_{\BOmega}})^T,\;\bphi=\BQ\bpsi=\sum_{j=1}^{N_{\bv}}\omega_j\bpsi_j,\quad j=1,\dots,N_{\BOmega}.
\end{align}
\end{subequations}
Here, $\BQ=\sum_{j=1}^{N_{\BOmega}}\omega_j\bpsi_j\in\mathbb{R}^{N_{\BOmega}N_{\textrm{DOF}}\times N_{\textrm{DOF}}}$ is the discrete integration operator.
The discrete advection operator $\BD_j\in\mathbb{R}^{N_{\textrm{DOF}}\times N_{\textrm{DOF}}}$, the discrete scattering and total cross sections $\BSigma_s,\BSigma_t\in\mathbb{R}^{N_{\textrm{DOF}}\times N_{\textrm{DOF}}}$, the discrete source and boundary flux $\bG,\bg_j^{(\textrm{bc})}\in\mathbb{R}^{N_\textrm{DOF}}$ are defined as:
\begin{subequations}
\begin{align}
&(\BD_{j})_{kl} = -\sum_{i=1}^{N_x}\int_{\mathcal{T}_i} (\BOmega_j\cdot\nabla\eta_k(\bx)) \eta_l(\bx) d\bx+ \sum_{\mathcal{E}\in\partial\mathcal{T}_h\setminus\partial\mathcal{T}_{h,j}^{(\textrm{ibc})}}\int_{\mathcal{E}} \widehat{\BH}\left(\BOmega_j, \eta_l,\bn(\bx)\right)\eta_k(\bx) d\bx,
\\
&(\BSigma_t)_{kl} = \sum_{i=1}^{N_x}\int_{\mathcal{T}_i} \sigma_t(\bx)\eta_k(\bx) \eta_l(\bx) d\bx,
\quad (\BSigma_s)_{kl} = \sum_{i=1}^{N_x}\int_{\mathcal{T}_i} \sigma_s(\bx)\eta_k(\bx) \eta_l(\bx) d\bx,
\\
&(\bG)_k= \sum_{i=1}^{N_x}\int_{\mathcal{T}_i} G(\bx)\pzc{\eta}_k(\bx) d\bx,\quad
(\bg_j^{(\textrm{bc})})_k = -\sum_{\mathcal{E}\in\partial\mathcal{T}_{h,j}^{(\textrm{ibc})}}\int_{\mathcal{E}} g(\bx)\pzc{\eta_k}(\bx) \BOmega_j\cdot \bn(\bx) d\bx,\\
&\wbG_j=\bG+\bg_j^{\textrm{\pzc{(bc)}}}.
\end{align}
\end{subequations}
Equation \eqref{eq:dg_matrix_vec} can be further rewritten \pzc{as} a fully coupled system:
\begin{equation}
\label{eq:one_equation}
\BA \bpsi =
\left(\begin{matrix}
\BD_1+\BSigma_t - \omega_1\BSigma_s & -\omega_2 \BSigma_s & \dots & -\omega_{N_{\BOmega}}\BSigma_s\\
-\omega_1\BSigma_s & \BD_2+\BSigma_t-\omega_2\BSigma_s & \dots & -\omega_{N_{\BOmega}}\BSigma_s \\
\vdots & \vdots & \vdots & \vdots \\
-\omega_1\BSigma_s & -\omega_2\BSigma_s & \dots & \BD_{N_{\BOmega}}+\BSigma_t-\omega_{N_{\BOmega}}\BSigma_s
\end{matrix}\right)
\left(
\begin{matrix}
\bpsi_1\\
\bpsi_2\\
\vdots\\
\bpsi_{N_{\BOmega}}
\end{matrix}
\right)
=
\left(\begin{matrix}
\wbG_1\\
\wbG_2\\
\vdots\\
\wbG_{N_{\BOmega}}
\end{matrix}
\right)=\bb.
\end{equation}

%%%%%%%%%%%%%%%%%%%%%%%%%%%%%%%%%%%%%%%%%%%%%%%%%%%%%%%%%%%%
\section{Krylov Solver and Source Iteration with Classical Synthetic Acceleration\label{sec:sa-review}}
In this section, we first briefly review the basic idea of Synthetic Acceleration (SA) under the Source Iteration (SI) framework \cite{adams2001discontinuous,azmy2010advances}. Then, following \cite{warsa2004krylov}, we review how to extend SA to Krylov methods for  faster convergence and enhanced robustness.

%%%%%%%%%%%%%%%%%%%%%%%%%%%%%%%%%%%
\subsection{Source Iteration and Synthetic Acceleration}
Each iteration of SI-SA comprises \pzc{an} SI step and \pzc{an} SA step. 
In the $l$-th iteration, to avoid directly solving the fully coupled system \eqref{eq:one_equation}, \textbf{the SI step} updates the angular flux $\psi$ by freezing the scalar flux $\phi$:
\begin{equation}
(\BD_j+\BSigma_t)\bpsi_{j}^{(l)} = \BSigma_s \bphi^{(l-1)}+\wbG_j,
\quad \bphi^{(l,*)}=\sum_{j=1}^{N_{\BOmega}}\omega_j \bpsi^{(l)}_j,\quad \bphi^{(0)}\text{ is given by an initial guess}.
\label{eq:SI-equation}
\end{equation} 
With proper choices of basis functions in the upwind DG spatial discretization, the linear system $\BD_j+\BSigma_s$ is block lower triangular if all elements are reordered along the upwind direction for $\BOmega_j$. After this reordering, equation \eqref{eq:SI-equation} for each $\BOmega_j$ can be solved with one block Gauss-Siedel iteration that sweeps along the upwind direction, namely transport sweep \cite{Adams2002FastIM}.  In practice, transport sweep is implemented in a matrix-free manner. 
Utilizing Fourier analysis, SI is proved to  converge, but the convergence speed may be arbitrarily slow for scattering dominant problems \cite{Adams2002FastIM}.

\textbf{The SA step} accelerates the convergence of SI by introducing a correction to the scalar flux after each source iteration
$\bphi^{(l)}=\bphi^{(l,*)}+\delta\bphi^{(l)}$. 

\pzc{To derive the ideal scalar flux correction, we first introduce the ideal angular flux correction $\delta \bpsi_j^{(l)}=\bpsi_j-\bpsi_j^{(l)}$, which is the difference between the true solution, $\bpsi_j$, and the  solution at the $l$-th iteration, $\psi_{j}^{(l)}$.} \pzc{By subtracting equation \eqref{eq:SI-equation} from equation \eqref{eq:discrete-rte}, we show that $\delta\bpsi_j^{(l)}$} satisfies an ideal discrete kinetic correction equation:
\begin{equation}
\label{eq:discrete_correction}
 (\BD_{j}+\BSigma_t) \delta\bpsi_{j}^{(l)} = \BSigma_s \delta\bphi^{(l)}+ \BSigma_s\left(\bphi^{(l,*)}-\bphi^{(l-1)}\right),\quad    \delta\bphi^{(l)}=\sum_{j=1}^{\pzc{N_{\BOmega}}}\omega_j \delta\bpsi_j^{(l)}.
\end{equation}
\pzc{Here, one can show that $\delta\bphi^{(l)}=\sum_{j=1}^{N_{\BOmega}}\omega_j\delta\bpsi_j^{(l)}$ is the ideal scalar flux correction satisfying $\delta\bphi^{(l)}=\bphi-\bphi^{(l,*)}$.} 
This ideal correction equation can be rewritten as 
\begin{equation}
\BA\delta\bpsi^{(l)} = \delta \bb^{(l)}=\underbrace{\left(\big(\BSigma_s(\bphi^{(l,*)}-\bphi^{(l-1)})\big)^T,\dots,\big(\BSigma_s(\bphi^{(l,*)}-\bphi^{(l-1)})\big)^T\right)^T}_{\text{Repeat } N_{\BOmega}\text{\;times}}.\label{eq:one_equation_correction}
\end{equation}
Using this ideal correction, SI converges in the next iteration. However, solving the ideal kinetic correction equation  \eqref{eq:discrete_correction} is equally expensive as solving the original linear system \eqref{eq:one_equation}. 

In practice, instead of solving the ideal kinetic correction equation \eqref{eq:discrete_correction}, we solve a computationally cheap low order approximation to it: 
\begin{equation}
\BC\delta\bphi^{(l)}=\BSigma_s(\bphi^{(l,*)}-\bphi^{(l-1)}),\label{eq:low-order-correction}
\end{equation}
and then correct the scalar flux as
\begin{equation}
\bphi^{(l)}=\bphi^{(l,*)}+\delta\bphi^{(l)}=\bphi^{(l,*)}+\BC^{-1}\BSigma_s(\bphi^{(l,*)}-\bphi^{(l-1)}).
\label{eq:practical-SA}
\end{equation}
The key of SA is to find a good low order approximation \eqref{eq:low-order-correction}.
Diffusion Synthetic Acceleration (DSA) \cite{kopp1963synthetic,alcouffe1977dittusion,adams1992diffusion,wareing1993new,Adams2002FastIM} exploits the diffusion limit of the ideal correction equation and a discretization partially/fully consistent with the upwind DG discretization (see \ref{sec:dsa} for details).  
Quasi-Diffusion methods \cite{gol1964quasi,anistratov1993nonlinear,olivier2023family} use the variable Eddington factor, while S$2$SA \cite{lorence1989s} utilizes the $S_2$ angular discretization. 

\textbf{Memory-efficient SI-SA:} Equation \eqref{eq:SI-equation} gives 
\begin{equation}
    \bpsi_j^{(l)} = (\BD_j+\BSigma_t)^{-1}\BSigma_s\bphi^{(l-1)}+(\BD_j+\BSigma_t)^{-1}\wbG_j.
\end{equation}
Numerically integrating equation \eqref{eq:SI-equation} in the angular space, we obtain 
\begin{equation}
    \bphi^{(l,*)}=\BT\BSigma_s\bphi^{(l-1)}+\wbb,\quad \BT=\sum_{j=1}^{N_{\BOmega}}\omega_j(\BD_j+\BSigma_t)^{-1},\quad\wbb=\sum_{j=1}^{N_{\BOmega}}\omega_j(\BD_j+\BSigma_t)^{-1}\wbG_j.\label{eq:memory-efficient-SI}
\end{equation}
Here, the operation of the operator $\BT=\sum_{j=1}^{N_{\BOmega}}\omega_j(\BD_j+\BSigma_t)^{-1}$ is implemented based on numerical integration and matrix-free transport sweeps.
Utilizing equation \eqref{eq:memory-efficient-SI}, SI-SA can be rewritten into  a memory-efficient formulation only involving the scalar flux:
\begin{enumerate}
    \item SI update: $\bphi^{(l,*)}=\BT\BSigma_s\bphi^{(l-1)}+\wbb$;
    \item SA correction: $\bphi^{(l)}=\bphi^{(l,*)}+C^{-1}\BSigma_s(\bphi^{(l,*)}-\bphi^{(l-1)})$.
\end{enumerate}
With this formulation, there is no need to explicitly save the angular flux $\bpsi_j$, however, the computational cost is almost the same as the original SI-SA in \eqref{eq:SI-equation} due to the numerical integration and transport sweeps in the discrete operator $\BT=\sum_{j=1}^{N_{\BOmega}}\omega_j(\BD_j+\BSigma_t)^{-1}$.

%%%%%%%%%%%%%%%%%%%%%%%%%%%%%%%%%%%
\subsection{GMRES with SA based preconditioner\label{sec:sa-gmres}}
\textbf{SI-SA as a Richardson iteration with left preconditioner:}
Following \cite{warsa2004krylov},  a preconditioner for discrete RTE can be constructed based on SA. Similar to the derivation of equation \eqref{eq:memory-efficient-SI}, we can derive the following equation from \eqref{eq:dg_matrix_vec} :
\begin{equation}
   \WBA\bphi=(\BI-\BT\BSigma_s)\bphi=\big(\BI-\big(\sum_{j=1}^{N_{\BOmega}}\omega_j(\BD_j+\BSigma_t)^{-1}\big)\BSigma_s\big)\bphi =\wbb\quad\text{and} \quad\wbb=\sum_{j=1}^{N_{\BOmega}}\omega_j(\BD_j+\BSigma_t)^{-1}\wbG_j.
\label{eq:memory-efficient-rte}
\end{equation}
\pzc{classical} SI without SA correction can be written as a Richardson iteration for \eqref{eq:memory-efficient-rte}:
\begin{align}
\bphi^{(l)}=\bphi^{(l,*)}&=\BT\BSigma_s\bphi^{(l-1)}+\wbb\notag\\
&= \bphi^{(l-1)}+(\wbb-(\BI-\BT\BSigma_s)\bphi^{(l-1)})\notag\\
&=\bphi^{(l-1)}+\br^{(l-1)}.
\label{eq:richardson}
\end{align}
Equation \eqref{eq:richardson} implies that the residual of the $(l-1)$-th iteration of the memory-efficient SI, $\br^{(l-1)}$, satisfies
\begin{equation}
\br^{(l-1)}=\bphi^{(l,*)}-\bphi^{(l-1)}.
    \label{eq:SI-residual-relation}
\end{equation}
Utilizing equation \eqref{eq:practical-SA}, \eqref{eq:memory-efficient-SI} and \eqref{eq:SI-residual-relation}, we rewrite SI-SA as
\begin{align}
\bphi^{(l)}&=\bphi^{(l,*)}+\BC^{-1}\BSigma_s(\bphi^{(l,*)}-\bphi^{(l-1)})\notag\\
&=\bphi^{(l,*)}+\BC^{-1}\BSigma_s\br^{(l-1)}\notag\\
&= \bphi^{(l-1)}+\br^{(l-1)}+\BC^{-1}\BSigma_s\br^{(l-1)}\notag\\
&=\bphi^{(l-1)}+(\BI+\BC^{-1}\BSigma_s)\br^{(l-1)}.
\label{eq:preconditioned-richardson}
\end{align}
Here, $(\BI+\BC^{-1}\BSigma_s)\br^{(l-1)}=(\BI+\BC^{-1}\BSigma_s)\left(\wbb-(\BI-\BT\BSigma_s)\bphi^{(l-1)}\right)$ is the residual of the left preconditioned  system 
\begin{equation}
(\BI+\BC^{-1}\BSigma_s)(\BI-\BT\BSigma_s)\bphi=(\BI+\BC^{-1}\BSigma_s)\wbb.
\end{equation}
Comparing equation \eqref{eq:richardson} and equation \eqref{eq:preconditioned-richardson}, we can see that introducing the SA step is equivalent to applying the left preconditioner $\BI+\BC^{-1}\BSigma_s$.

\textbf{SA Preconditioned GMRES:} 
As discussed in \cite{azmy2002unconditionally,warsa2004krylov}, SI-DSA may lose its efficiency and robustness for challenging problems with discontinuous materials.  
To accelerate the convergence and improve the robustness, SA preconditioners are extended to Krylov methods solving $(\BI-\BT\BSigma_s)\bphi=\wbb$. Inspired by viewing SI-SA as a left preconditioned Richardson iteration, the left or right preconditioner for Krylov methods can be chosen as  $\BI+\BC^{-1}\BSigma_s$, where $\BC$ is the low order approximation to the ideal correction equation in the underlying SA.
%%%%%%%%%%%%%%%%%%%%%%%%%%%%%%%%%%
\subsection{Limitations of DSA}
Classical DSA approximates the ideal kinetic correction equation \eqref{eq:discrete_correction} with its diffusion limit. However, it may become less effective, when the problem is far from its diffusion limit \cite{ren2019fast}. Moreover, the low rank structures of the solution manifold for the parametric problem is not leveraged.  
%%%%%%%%%%%%%%%%%%%%%%%%%%%%%%%%%%%%%%%%%%%%%%%%%%%%%%%%%%%%
\section{Reduced order model enhanced preconditioner for Krylov method\label{sec:romsad-preconditioner}}
To address the limitations of the classical DSA preconditioner, we exploit a ROM for the ideal correction equation to design a new SA preconditioner. This ROM is based on the original kinetic description of the ideal correction equation and leverages low rank structures in parametric problems.
It is constructed following the offline-online decomposition approach. In the offline stage, we build the ROM by extracting low rank structures from data, specifically solutions to the ideal correction equation for parameters in a training set. In the online stage, to obtain high-fidelity solutions for a new parameter outside the training set, we apply a  Krylov solver accelerated by the constructed ROM-enhanced preconditioner. Specifically, following our previous work to accelerate SI \cite{peng2024romsad}, we will use a  hybrid SA preconditioner that utilizes ROM-based corrections in the first few iterations to improve efficiency and then switches to DSA for better robustness.  Due to this switching, the resulting preconditioner is nonlinear. Hence, we choose flexible general minimal residual method (FGMRES) \cite{saad1986gmres}, allowing the use of nonlinear preconditioners, as our underlying Krylov solver. 

In this section, we \pzc{introduce affine parametric problems and outline basic ideas of ROMs for such problems in Sec. \ref{sec:parametric-assumption}}, review the algorithm details of FGMRES in Sec. \ref{sec:fgmres}, \pzc{define} our ROM-enhanced preconditioner in Sec. \ref{sec:romsad},  and discuss how to construct the underlying ROM  in Sec. \ref{sec:rom-construction}. 
%%%%%%%%%%%%%%%%%%%%%%%%%%%%%%%%%%%
\subsection{\pzc{ROM for affine parametric problem \label{sec:parametric-assumption}}}
\pzc{\textbf{Affine parametric problem.} Here, we introduce our notations and assumptions. Let the full order problem \eqref{eq:one_equation} for the parameter $\bmu$ be $\BA_{\bmu}\bpsi_{\bmu}=\bb_{\bmu}$, and the equivalent memory-efficient formulation be 
$\WBA_{\bmu}\bphi_{\bmu}=\wbb_{\bmu}$. Here, $\bmu$ is the parameter of the parametric problem, such as the parametrization of material properties or boundary conditions. The subscript $\bmu$ denotes the $\bmu$-dependence. Additionally, for simplicity, we assume the operator $\BA_\bmu$ depends on the parameter $\bmu$ in an affine manner, i.e. 
\begin{equation}
\BA_{\bmu}=\sum_{s=1}^{m}a_s(\bmu)\BA_s,  
\label{eq:affine-assumption}
\end{equation} 
where $\BA_s$'s are constant matrices independent of $\bmu$.} 

\pzc{\textbf{ROM construction.} We follow an offline-online decomposition framework to construct the ROM. In the offline stage, we find low rank structures in the ideal correction equation for the kinetic parametric problems, i.e. $\BA_{\bmu}\delta\bpsi_{\bmu}^{(l)}=\delta\bb^{(l)}_{\bmu}$, from data. Exploiting these low-rank structures, we construct a reduced order space whose bases are column vectors of $\BU_r\in\mathbb{R}^{N_{\BOmega}N_{\textrm{DOF}}\times r}$. In the online stage, a reduced order solution to the ideal kinetic correction equation \eqref{eq:one_equation_correction} can be efficiently obtained by  projecting  this equation onto the reduced order space $\BU_r$:
\begin{equation}
\delta\bpsi_{\bmu}^{(l)}\approx \BU_r \bc_{\bmu}^{(l)}, \quad \BU_r^T\BA_{\bmu}\BU_r \bc_{\bmu,r}^{(l)}=\BU_r^T\delta\bb_{\bmu}^{(l)}.
\end{equation}
}
\pzc{Under the affine assumption, the reduced order operator for any $\bmu$ can be fast constructed online as 
\begin{equation}
\BU_r^T\BA_{\bmu}\BU_r=\sum_{s=1}^ma_s(\bmu)\BA_{s,r},
\end{equation} where $\BA_{s,r}=\BU_r^T\BA_s\BU_r$ are precomputed and saved offline.}

\pzc{A key point in our method is that the ROM is directly built for the ideal kinetic correction equation $\BA_{\bmu}\delta\bpsi_{\bmu}=\delta\bb$, rather than its memory-efficient formulation, $(\BI-\BT\BSigma_s)\delta\bphi^{(l)}=\delta\widetilde{\bb}$. Projecting $\BA_{\bmu}$ avoids transport sweeps required when projecting $\BI-\BT\BSigma_s$. Moreover, the operator $\BT=\sum_{j=1}^{N_{\BOmega}}\omega_j(\BD_j+\BSigma_t)^{-1}$ exhibits a non-affine dependence on the total cross section $\sigma_t$, whenever  $\sigma_t$ is parametric.  This non-affine dependence can lead to a loss of online efficiency, as the projected reduced-order operator cannot be rapidly constructed online using affine decomposition and offline precomputations \cite{behne2022minimally}.  To address efficiency losses due to non-affine dependence, classical approaches such as the empirical interpolation method (EIM) \cite{barrault2004empirical} and discrete EIM (DEIM) \cite{chaturantabut2010nonlinear} rely on sub-sampling rows of reduced-order matrices and vectors. However, when the matrix-vector product determined by $\BT$ is implemented via a matrix-free transport sweep, computing the outputs of the sampled rows requires first computing all ``previous rows along the upwind direction". Consequently, when the matrix-free transport sweep is applied, simply applying EIM and DEIM may not achieve the desired online efficiency.}

\begin{rem}
\pzc{\textbf{Extensions to non-affine problems.} The focus of this paper is to demonstrate the potential of using ROMs to design more efficient preconditioners for Krylov methods solving parametric RTE. For simplicity, we only consider affine problems. However, our current strategy directly building ROMs on the kinetic level can be combined with classical hyper-reduction techniques, such as EIM \cite{barrault2004empirical} and DEIM \cite{chaturantabut2010nonlinear}, to build efficient ROMs for non-affine  and nonlinear problems.}
\end{rem}
%%%%%%%%%%%%%%%%%%%%%%%%%%%%%%%%%%%
\subsection{Flexible GMRES\label{sec:fgmres}}
Following \cite{saad1993flexible}, we outline FGMRES in Alg. \ref{alg:fgmres} which allows using a nonlinear preconditioner $\BM_j^{-1}$. Here, the subscript $j$ indicates the iteration dependence of the nonlinear preconditioner.  The subtle but crucial difference between FGMRES and right preconditioned GMRES is that the preconditioned vectors $\bz^{(j)}$ are saved and utilized to update the solution. If the preconditioner  is linear, i.e. $\BM_j^{-1}=\BM^{-1}$ for any $j$, FGMRES is mathematically equivalent to right preconditioned GMRES.  
%%%%%%%%%%%%%%%%%%%%%%%%%%%
\begin{algorithm}[H]
\caption{Flexible GMRES method to solve $\WBA\bphi=\wbb$ with a  nonlinear preconditioner $\BM_j^{-1}$ \cite{saad1993flexible}. \label{alg:fgmres} }
\begin{algorithmic}[1]
\STATE{Given an initial guess $\bphi^{(0)}$, a nonlinear preconditioner $\BM_j^{-1}$, and the maximum number of iterations allowed $m$.}
\STATE{\textit{\textbf{Initialization:}}} 
\IF{$||\bphi^{(0)}||_2==0$}
    \STATE{$\br^{(0)}=\wbb$,}
    \ELSE
    \STATE{$\br^{(0)}=\wbb-\WBA\bphi^{(0)}$.}
\ENDIF
\STATE{Define $\beta=||\br^{(0)}||_2$ and $\bv^{(1)}=\br^{(0)}/\beta$. Allocate the memory for the Hessenberg matrix $\BH\in\mathbb{R}^{(m+1)\times m}$.}
\STATE{\textit{\textbf{Anorldi process:}}}\FOR{ $j=1:m$}
\STATE{\textbf{Apply the nonlinear preconditioner:}  $\bz^{(j)}:=\BM_j^{-1}\bv^{(j)}$.}
\STATE{ Compute $\bw:=\WBA\bz^{(j)}$.}
    \FOR{$i=1:j$}
        \STATE{$\BH_{ij}=\bw^T \bv^{(i)}$, and update 
        $\bw=\bw-\BH_{ij}\bv^{(i)}$.}
    \ENDFOR
\STATE{Compute $\BH_{j+1,j}=||\bw||_2$.}
\STATE{Update $\bv^{(j+1)}=\bw/\BH_{j+1,j}$.}
\IF{Stopping criteria satisfied}
    \STATE{Stop iterations.}
\ENDIF
\ENDFOR
\STATE{\textit{\textbf{Obtain the solution:}}} 
\STATE{Define $\widetilde{\BZ}=\left(\bz^{(1)},\bz^{(2)},\dots,\bz^{(j)}\right)$ and $\widetilde{\BH}=\left(\BH_1,\BH_2,\dots,\BH_j\right)$. Solve the minimization problem $\by^{*}=\arg\min_{\by\in\mathbb{R}^j}||\beta\be_1-\widetilde{\BH}\by||_2$, where $\be_1=(1,0,\dots,0)^T\in\mathbb{R}^j$.}
\STATE{Compute and return $\bphi:=\bphi^{(0)}+\widetilde{\BZ}\by^*$.}
\end{algorithmic}
\end{algorithm}

%%%%%%%%%%%%%%%%%%%%%%%%%%%%%%%%%%%
\subsection{ROM-enhanced SA preconditioner \label{sec:romsad}}
To leverage the original kinetic description of the correction equation and low rank structures in parametric problems, we propose a SA strategy combining ROM-based SA and DSA, called ROMSAD, under the SI framework in \cite{peng2024romsad}. Here, we briefly review its basic ideas and outline how to extend it to the Krylov framework.

\pzc{Let columns vectors of $\BU_r\in\mathbb{R}^{N_{\BOmega}N_{\textrm{DOF}}\times r}$ be an orthornormal basis of our reduced order space.} We denote rows corresponding to degrees of freedom for $\psi_h(\cdot,\Omega_j)$ in $\BU_r=(\BU_{r,1}^T,\dots,\BU_{r,N_{\BOmega}}^T)^T$ as $\BU_{r,j}\in\mathbb{R}^{N_{\textrm{DOF}}\times r}$. In the online stage, we seek the correction by projecting the ideal kinetic correction equation \eqref{eq:one_equation_correction} onto the reduced order space $\BU_r$:
\begin{subequations}
  \label{eq:rom-correction}
\begin{align}
&\BU_r^T\BA\BU_r\bc^{(l)}_r = \BU_r^T \delta\bb^{(l)}=\sum_{j=1}^{N_{\BOmega}}\BU_{r,j}^T\BSigma_s(\bphi^{(l,*)}-\bphi^{(l-1)})=(\sum_{j=1}^{N_{\BOmega}}\BU_{r,j})^T\BSigma_s\br^{(l-1)},\quad\bc_r^{(l)}\in\mathbb{R}^{r},\\
&\delta \bpsi^{(l)}\approx \BU_r\bc^{(l)}_r,\quad \delta\bphi^{(l)}=\sum_{j=1}^{N_{\BOmega}}\omega_j\delta\bpsi_j^{(l)}\approx\sum_{j=1}^{N_{\BOmega}}\omega_j \BU_{r,j}\bc^{(l)}_r=\big(\sum_{j=1}^{N_{\BOmega}}\omega_j\BU_{r,j}\big)\bc_r^{(l)}.
\end{align}
\end{subequations}
Our ROMSAD method uses ROM-based corrections in \eqref{eq:rom-correction} for the first few iterations and then switches to DSA. We need this switching, as the ROM is constructed only using the data for the first few iterations, i.e. $\delta\bpsi^{(l)}$, $1\leq l\leq\mfw$,  to avoid excessive memory costs due to the high dimensionality of the angular flux. 

To extend ROMSAD method to the Krylov framework, we need to derive the correction operator $``\BC^{-1}"$ in \eqref{eq:practical-SA} for the ROM-based correction. Based on  \eqref{eq:rom-correction}, we have 
\begin{align}
\delta\bphi^{(l)}=\sum_{j=1}^{N_{\BOmega}}\omega_j\delta\bpsi^{(l)}_j&\approx 
(\sum_{j=1}^{N_{\BOmega}}\omega_j \BU_{r,j})(\BU_r^T\BA\BU_r)^{-1}\big(\sum_{j=1}^{N_{\BOmega}}\BU_{r,j})^T\BSigma_s\br^{(l-1)}.
    \label{eq:C-rom-SA}
\end{align}
Comparing \eqref{eq:practical-SA} and \eqref{eq:C-rom-SA}, we conclude  $\BC_{\textrm{ROM}}^{-1}= \left(\sum_{j=1}^{N_{\BOmega}}\omega_j \BU_{r,j}\right)(\BU_r^T\BA\BU_r)^{-1}\left(\sum_{j=1}^{N_{\BOmega}}\BU_{r,j}\right)^T$. Following Sec. \ref{sec:sa-gmres},  the SA preconditioner corresponding to this ROM correction is
\begin{equation}
\BM^{-1}_{\textrm{ROM}}=\BI+\BC_{\textrm{ROM}}^{-1}\BSigma_s = \BI+(\sum_{j=1}^{N_{\BOmega}}\omega_j \BU_{r,j})(\BU_r^T\BA\BU_r)^{-1}\big(\sum_{j=1}^{N_{\BOmega}}\BU_{r,j}\big)^T\BSigma_s.\label{eq:rom-sa-preconditioner}
\end{equation}
Then, the nonlinear ROMSAD preconditioner for the $j$-th iteration of FGMRES can be defined as
\begin{equation}
    \BM^{-1}_j = \begin{cases}
                \BM_{\textrm{ROM}}^{-1},\quad 1\leq j\leq \mfw,\\
                \BM_{\textrm{DSA}}^{-1},\quad\text{otherwise}.
                \end{cases}
                \label{eq:ROMSAD-preconditioner}
\end{equation}  
Here, the window size, $\mfw$, is defined as the number of iterations whose data are utilized in the ROM construction, and $\BM_{\textrm{DSA}}^{-1}$ is the DSA-preconditioner. We want to point out that a more complicated switching strategy is used in our previous work for SI \cite{peng2024romsad}. Due to the better robustness of Krylov method, the switching strategy in \eqref{eq:ROMSAD-preconditioner} works well for all our numerical tests.

%%%%%%%%%%%%%%%%%%%%%%%%%%%%%%%%%%%

%%%%%%%%%%%%%%%%%%%%%%%%%%%%%%%%%%%
\begin{rem}
Our ROM-based SA preconditioner is different from directly extending the ROM-based preconditioner for elliptic equations \cite{santo2018multi} to RTE, as the starting points for these two methods are different. Based on the SA framework, our ROM-based SA preconditioner is 
$$\BI+\sum_{j=1}^{N_{\BOmega}}(\omega_j \BU_{r,j})(\BU_r^T\BA\BU_r)^{-1}\big(\sum_{j=1}^{N_{\BOmega}}\BU_{r,j}\big)^T\BSigma_s.$$ \cite{santo2018multi} derives a ROM-based preconditioner by replacing the ``coarse grid space" in two-grid methods with a reduced order space. Directly applying \cite{santo2018multi} to the memory-efficient formulation of RTE \eqref{eq:memory-efficient-rte} results in the preconditioner $\widetilde{\BU}_r\left(\widetilde{\BU}_r^T(\BI-\BT\BSigma_t)\widetilde{\BU}_r\right)^{-1}\widetilde{\BU}_r^T$ where $\widetilde{\BU}_r$ is the reduced order basis for a ROM of the memory-efficient formulation. Clearly, this direct extension of \cite{santo2018multi} may suffer from  the efficiency loss due to the aforementioned non-affine issue \pzc{in Sec. \ref{sec:parametric-assumption}}.
\end{rem}

%%%%%%%%%%%%%%%%%%%%%%%%%%%%%%%%%%%
\subsection{ROM construction \label{sec:rom-construction}}
The only remaining question  is how to build the ROM for the ideal kinetic correction equation. 
To answer this question, we first figure out  the ideal kinetic correction equation for FGMRES and then rewrite it into an equivalent form to obtain its solution efficiently without solving it. Then, we build the ROM using a greedy algorithm \cite{hesthaven2016certified} that iteratively identifies representative parameters from a training set based on a residual-based error indicator. Using this strategy, we only need to generate the high-fidelity ideal angular flux corrections for these identified  parameters, while  our previous work \cite{peng2024romsad}, constructed the ROM by proper orthogonal decomposition (POD) using high-fidelity data for all parameters in a training set. Consequently, the underlying ROM may be constructed more efficiently compared to \cite{peng2024romsad}. 
%%%%%%%%%%%%%%%%%%%%%%%%%%%%%%%%%%%
\subsubsection{Ideal correction under the Krylov framework and data generation \label{sec:data-generation}}
As discussed in the last paragraph of Sec. \ref{sec:romsad}, we need to construct the ROM for the kinetic correction equation instead of its memory-efficient formulation
to avoid the non-affine issues. Here, we identify the ideal kinetic correction equation for FGMRES as follows.

Recall that, under the SI framework, the residual for the $(l-1)$-th iteration, $\br^{(l-1)}$, satisfies  $\br^{(l-1)}=\bphi^{(l,*)}-\bphi^{(l-1)}$. Substituting this relation into equation \eqref{eq:SI-equation},  the ideal angular flux  correction for SI can be rewritten as: 
\begin{equation}
(\BD_j+\BSigma_t)\delta\bpsi_j^{(l)}
=\BSigma_s\delta\bphi^{(l)}+\BSigma_s\br^{(l-1)},\quad \delta\bphi^{(l)}=\sum_{j=1}^{N_{\BOmega}}\omega_j\delta\bpsi_{j}^{(l)},\quad 1\leq j\leq N_{\BOmega}.
\end{equation}
In the preconditioning step of FGMRES, the residual in SI, $\br^{(l-1)}$, is  replaced by the orthogonal basis of the Krylov space $\bv^{(l)}$. Therefore, the ideal angular flux correction for FGMRES, namely $\delta\bpsi_j^{(l)}$, solves 
\begin{equation}
(\BD_j+\BSigma_t)\delta\bpsi_j^{(l)}=\BSigma_s\delta\bphi^{(l)}+\BSigma_s \bv^{(l)},\quad \delta\bphi^{(l)}=\sum_{j=1}^{N_{\BOmega}}\omega_j\delta\bpsi^{(l)}_j,\quad j=1,\dots,N_{\BOmega}.  
\label{eq:ideal-kinetic-correction-fgmres}
\end{equation}
\pzc{Solving equation \eqref{eq:ideal-kinetic-correction-fgmres} is as expensive as directly solving $\WBA\bphi=\wbb$. To build the ROM for equation \eqref{eq:ideal-kinetic-correction-fgmres} without solving it, we reformulate it as follows.}
\pzc{Define $\bxi^{(l)}$ as the solution to 
\begin{equation}
    \WBA\bxi^{(l)}=(\BI-\BT\BSigma_s)\bxi^{(l)}=\bv^{(l)}.\label{eq:xi-def}
\end{equation}
Substituting equation \eqref{eq:xi-def} into
\eqref{eq:ideal-kinetic-correction-fgmres}, we obtain
\begin{align}
    (\BD_j+\BSigma_t)\delta\bpsi_j^{(l)}&=
    \BSigma_s\delta\bphi^{(l)}+\BSigma_s(\BI-\BT\BSigma_s)\bxi^{(l)},\\
    \delta\bpsi_j^{(l)}&=(\BD_j+\BSigma_t)^{-1}\left(\BSigma_s\delta\bphi^{(l)}+\BSigma_s(\BI-\BT\BSigma_s)\bxi^{(l)}\right).
\end{align}
Then, the scalar flux correction, $\delta\bphi^{(l)}$, can be obtained through numerical quadrature as
\begin{align}
\delta\bphi^{(l)}&=\big(\sum_{j=1}^{N_{\BOmega}}\omega_j (\BD_j+\BSigma_t)^{-1}\big)\left(\BSigma_s\delta\bphi^{(l)}+\BSigma_s(\BI-\BT\BSigma_s)\bxi^{(l)}\right).
\end{align}
Based on the definition of $\BT=\sum_{j=1}^{N_{\BOmega}}\omega_j (\BD_j+\BSigma_t)^{-1}$ in \eqref{eq:memory-efficient-SI}, we obtain
\begin{subequations}
\begin{align}
    \delta\bphi^{(l)}&= \BT\BSigma_s\delta\bphi^{(l)}+\BT\BSigma_s(\BI-\BT\BSigma_s)\bxi^{(l)}\\
    (\BI-\BT\BSigma_s)\delta\bphi^{(l)}&=\BT\BSigma_s(\BI-\BT\BSigma_s)\bxi^{(l)}\\
    &=\left(\BT\BSigma_s-(\BT\BSigma_s)^2\right)\bxi^{(l)}=(\BI-\BT\BSigma_s)\BT\BSigma_s\bxi^{(l)}.
\end{align}
\end{subequations}
}
\pzc{
Hence, 
\begin{equation}
\delta\bphi^{(l)}=\BT\BSigma_s\bxi^{(l)}.
\label{eq:density-relation}
\end{equation}
Substituting \eqref{eq:density-relation} and \eqref{eq:xi-def} into \eqref{eq:ideal-kinetic-correction-fgmres}, we have
\begin{align}
    (\BD_j+\BSigma_t)\delta\bpsi_j^{(l)}=\BSigma_s\BT\BSigma_s\bxi^{(l)}+\BSigma_s(\BI-\BT\BSigma_s)\bxi^{(l)}=\BSigma_s\bxi^{(l)}.
\end{align}
}
\pzc{
Hence, the ideal  correction equation \eqref{eq:ideal-kinetic-correction-fgmres} is equivalent to
\begin{subequations}
\label{eq:fgmres-psi-correction-whole}
\begin{align}
&(\BD_j+\BSigma_t)\delta\bpsi_j^{(l)}=\BSigma_s\bxi^{(l)},\quad j=1,\dots,N_{\BOmega},
\label{eq:fgmres-psi-correction}\\
&\WBA\pzc{\bxi}^{(l)}=(\BI-\BT\BSigma_s)\bxi^{(l)}=\bv^{(l)}.
\end{align} 
\end{subequations}
} 

Utilizing \eqref{eq:fgmres-psi-correction-whole}, snapshots for the ideal angular flux corrections can be efficiently constructed as long as we can  efficiently compute the solution to $\pzc{\WBA\bxi^{(l)}=(\BI-\BT\BSigma_s)\bxi^{(l)}=\bv^{(l)}}$ without directly solving this equation or using matrix-vector multiplications determined by $\WBA$. Similar to \cite{santo2018multi}, this can be achieved by utilizing the definition of FGMRES method in Alg. \ref{alg:fgmres}.
In the first iteration, we have
\begin{align}
&\WBA\pzc{\bxi}^{(1)}=(\BI-\BT\BSigma_s)\pzc{\bxi}^{(1)}=\bv^{(1)}=\br^{(0)}/\beta=\br^{(0)}/||\br^{(0)}||_2=(\wbb-\WBA\bphi^{(0)})/||\br^{(0)}||,\notag\\
&\pzc{\bxi}^{(1)}=\WBA^{-1}(\wbb-\WBA\bphi^{(0)})/||\br^{(0)}||_2=
\WBA^{-1}(\WBA\bphi-\WBA\bphi^{(0)})/||\br^{(0)}||_2=(\bphi-\bphi^{(0)})/||\br^{(0)}||_2,\label{eq:correction-fgmres1}
\end{align}
where $\bphi$ is the solution to $\WBA\bphi=(\BI-\BT\BSigma_s)\bphi=\wbb$, which can be well approximated by the converged solution of FGMRES. As a result, to compute $\pzc{\bxi}^{(1)}$, we only need to save the converged solution $\bphi$ and the initial guess $\bphi^{(0)}$.
Similarly,  based on line $12$ to line $17$ of Alg. \ref{alg:fgmres}, we can show that  $\pzc{\bxi}^{(l)}$ for $l\geq 2$ satisfying
\begin{align}
\WBA\pzc{\bxi}^{(l)}&=\bv^{(l)}=
(\WBA \bz^{(l-1)}-\sum_{i=1}^{l-1}\BH_{i,l-1}\bv^{(i)})/\BH_{l,l-1},\quad l\geq 2,\notag\\
\pzc{\bxi}^{(l)}&=\WBA^{-1}(\WBA \bz^{(l-1)}-\sum_{i=1}^{l-1}\BH_{i,l-1}\bv^{(i)})/\BH_{l,l-1}
=\frac{1}{\BH_{l,l-1}}\left(\bz^{(l-1)}-\sum_{i=1}^{l-1}\BH_{i,l-1}\WBA^{-1}\bv^{(i)}\right)\notag\\
&=\frac{1}{\BH_{l,l-1}}\left(\bz^{(l-1)}-\sum_{i=1}^{l-1}\BH_{i,l-1}\pzc{\bxi}^{(i)}\right),\quad l\geq 2.\label{eq:correction-fgmres2}
\end{align}
Utilizing \eqref{eq:correction-fgmres1} and \eqref{eq:correction-fgmres2}, we can iteratively compute $\pzc{\bxi}^{(l)}$ for the first $\mfw$ iterations, by saving the converged solution, namely $\pzc{\bphi}$, elements in the matrix $\BH$, vectors $\bz^{(l)}$ and  $\bv^{(l)}$.

In summary, under FGMRES framework, to generate the ideal angular flux corrections $\delta\bpsi^{(l)}$ for the first few iterations $1\leq l\leq \mfw$, we first compute $\pzc{\bxi}^{(l)}$ iteratively using \eqref{eq:correction-fgmres1} and \eqref{eq:correction-fgmres2}. Then, we construct the ideal angular flux corrections $\delta\bpsi^{(l)}$ by applying $\mfw$ additional transport sweeps to solve \eqref{eq:fgmres-psi-correction}. Using this two-step strategy, we generate $\delta\bpsi^{(l)}$ while avoiding directly solving \eqref{eq:ideal-kinetic-correction-fgmres} or using matrix-vector multiplications determined by $\WBA$.

%%%%%%%%%%%%%%%%%%%%%%%%%%%%%%%%%%%
\subsubsection{Greedy algorithm to build ROM\label{sec:greedy}}
In the offline stage, we will use a greedy algorithm  to build a ROM for the ideal kinetic correction equation for FGMRES \eqref{eq:ideal-kinetic-correction-fgmres}. The author is aware that a systematic study of greedy algorithms to build ROMs for parametric RTE is ongoing \cite{matsuda2025rbmrte}. The greedy algorithm iteartively enriches the reduced order space by sampling the most representative parameters from the training set. These representative parameters are identified by a residual-based error indicator. High-fidelity linear solves are only needed to generate data for these sampled parameters.

We use FGMRES with DSA preconditioner, which is equivalent to right preconditioned GMRES, and zero initial guess to generate required high-fidelity data. The steps of the greedy algorithm are as follows.
\begin{enumerate}
\item[\textbf{Input:}] The training set of  parameters $\mathcal{P}_{\textrm{train}}=\{\bmu_s\}_{s=1}^{N_{\textrm{train}}}$, maximum number of greedy iterations $N_{\textrm{g-iter}}$, the window size $\mfw$ determining high-fidelity data for how many iterations will be used, the initial sample $\bmu_1$ and the tolerance for the residual $\epsrom$,
\item \textbf{Preparation and initialization.} 
\begin{itemize}
    \item[] Use \eqref{eq:memory-efficient-SI} to compute the right-hand side of the memory-efficient formulation $\wbb_{\bmu}$ for all $\bmu\in\mathcal{P}_{N_{\textrm{train}}}$. Since zero initial guesses are applied, the initial residuals for FGMRES are $\br_{\bmu}^{(0)}=\wbb_{\bmu}=||\br_{\bmu}^{(0)}||_2\bv^{(1)}$. Compute and save $\boldeta_{\bmu}^{(0)}=\BSigma_s\br_{\bmu}^{(0)}$.
    Set  $k=1$, and the set of sampled parameters as $\mathcal{P}_{sampled}=\{\bmu_1\}$.
\end{itemize}
\item \textbf{Update reduced order basis.}
\begin{itemize}
    \item[(a)] \textbf{High fidelity solve for sampled parameters}. Use FGMRES to solve $\WBA_{\bmu_k}\bphi_{\bmu_k}=\wbb_{\bmu_k}$.  Use \eqref{eq:correction-fgmres1} and \eqref{eq:correction-fgmres2} to compute the solution to $\WBA\pzc{\bxi}^{(l)}_{\bmu_k}=\bv^{(l)}_{\bmu_k}$, $\pzc{\bxi}^{(l)}_{\bmu_k}$, for $1\leq l\leq \mfw$.
    \item[(b)] \textbf{Compute angular flux corrections}. Generate the angular flux correction $\delta\bpsi^{(l)}_{\bmu_k}$ for $1\leq l\leq \mfw$ by solving \eqref{eq:fgmres-psi-correction} with transport sweeps. 
    \item[(c)] \textbf{Update reduced order basis and operators.} Apply the modified Gram-Schmidt procedure \pzc{with a truncation step (see \ref{sec:mgs-procedure} for detials)}  to update $\BU_r$ by adding columns orthogonal to  its existing columns from the column space of the snapshot matrix $\BS_{\bmu_k}=[\delta\bpsi^{(1)}_{\bmu_k},\dots,\delta\bpsi^{(\mfw)}_{\bmu_k}]$, and then update reduced order operators correspondingly. If $k$ equal to $N_{\textrm{g-iter}}$, stop. Otherwise, go to step 3.
    \end{itemize}
\item \textbf{Greedy sampling.} 
    \begin{itemize}
    \item[(a)] \textbf{Generate reduced order solutions.}
    For all unsampled parameters, $\bmu\in\mathcal{P}_{\textrm{train}}\setminus\mathcal{P}_{\textrm{sampled}}$, compute the reduced order angular flux of the correction equation for the first iteration of FGMRES:
    \begin{subequations}
    \label{eq:rom-sampling}
    \begin{align}
    &\BU_r^T\BA_{\bmu}\BU_r\bc_{r,\bmu} =(\sum_{j=1}^{N_{\BOmega}}\BU_{r,j})^T\BSigma_s\br^{(0)}_{\bmu}=||\br^{(0)}_{\bmu}||_2(\sum_{j=1}^{N_{\BOmega}}\BU_{r,j})^T\BSigma_s\bv^{(1)}=(\sum_{j=1}^{N_{\BOmega}}\BU_{r,j}^T)\boldeta^{(0)}_{\bmu},\quad\bc_{r,\bmu}\in\mathbb{R}^{r}\label{eq:greedy-rom}\\
    &\delta\bpsi_{\bmu,\textrm{ROM}}\approx \BU_r\bc_r,\quad \delta\bphi_{\bmu,\textrm{ROM}}\approx\sum_{j=1}^{N_{\BOmega}}\omega_j \BU_{r,j}\bc_r. 
    \end{align}
    \end{subequations}
    \item[(b)] \textbf{Greedy sampling based on residuals.} Compute the residual of the current reduced order solutions
    \begin{align}
    \mathcal{R}_{\bmu}&=\max_{1\leq j\leq N_{\BOmega}}||(\BD_j+\BSigma_{t,\bmu})\delta\bpsi_{j,\textrm{ROM}}-\BSigma_{s,\bmu}\delta\bphi_{\textrm{ROM}}-\BSigma_{s,\bmu}\br_{\bmu}^{(0)}||_2\notag\\
    &=\max_{1\leq j\leq N_{\BOmega}}||(\BD_j+\BSigma_{t,\bmu})\delta\bpsi_{j,\textrm{ROM}}-\BSigma_{s,\bmu}\delta\bphi_{\textrm{ROM}}-\boldeta_{\bmu}^{(0)}||_2.\label{eq:rom-residual}
    \end{align}
    If $\max_{\bmu\in\mathcal{P}_{\textrm{train}}\setminus\mathcal{P}_{\textrm{sampled}}}\{\mathcal{R}_{\bmu}\}<\epsilon_{rom}$, stop the algorithm.
    Otherwise, sample $\bmu_{k+1}$ as  
$$\bmu_{k+1}={\arg\max}_{\bmu\in\mathcal{P}_{\textrm{train}}\setminus\mathcal{P}_{\textrm{sampled}}}\{\mathcal{R}_{\bmu}\},$$
     then update $\mathcal{P}_{\textrm{sampled}}:=\mathcal{P}_{\textrm{sampled}}\bigcup\{\bmu_{k+1}\}$ and $k:=k+1$. Go back to step 2. 
    \end{itemize}
\item{\textbf{Output:}} the reduced order basis $\BU_r$ \pzc{with $r=k\mfw$} and corresponding reduced order operators.
\end{enumerate}
In our greedy algorithm, representative parameters are sampled based on the correction equation only for the first iteration of FGMRES,  because given the initial guess and the right-hand side, they can be easily computed in the preparation step. 
Though we use zero initial guesses  in this paper, the greedy method can be straightforwardly extended to other initial guesses.

Compared with the POD method used in our previous work \cite{peng2024romsad}, the greedy algorithm reduces the number of high-fidelity solves with the additional cost of greedy sampling. Our numerical tests demonstrate that, when the number of representative parameters is smaller than the size of the training set, this additional costs is marginal compared to generating the high-fidelity data for the entire training set. 

\begin{rem}
\pzc{To some extent, the greedy algorithm can be seen as an error-indicator-pivoted QR decomposition. It iteratively constructs a set of truncated orthogonal basis, $\BU_r$, by only sampling as many parameters as necessary. When expanding the basis, it avoids adding less informative parameters corresponding to snapshots nearly linearly dependent on $\BU_r$. This is because the reduced order solution for such parameters lead to small residuals while the greedy algorithm samples the parameter with the largest residual.} 
\end{rem}
%%%%%%%%%%%%%%%%%%%%%%%%%%%%%%%%%%%%%%%%%%%%%%%%%%%%%%%%%%%%

%%%%%%%%%%%%%%%%%%%%%%%%%%%%%%%%%%%%%%%%%%%%%%%%%%%%%%%%%%%%
\section{Numerical results\label{sec:numerical}}
We will utilize a series of numerical examples in the 2D X-Y geometry to test the performance of FGMRES with the ROMSAD preconditioner (FGMRES-ROMSAD). We will compare its performance to GMRES with right DSA preconditioner (GMRES-DSA), since FGMRES with a linear preconditioner is mathematically equivalent to right preconditioned GMRES. \pzc{As discussed in \cite{saad2003iterative}, the main difference between left and right preconditioned GMRES is that the left preconditioned version minimizes the residual of the preconditioned system, while the right preconditioned version directly minimizes the original system. In practice, their performance is comparable.} Throughout this section, we use a linear DG space ($K=1$) for the spatial discretization. For the DSA preconditioner, we will use the fully consistent version described in \ref{sec:dsa}, unless otherwise specified. We observe that, in comparison with the partially consistent version, the fully consistent one results in less number of iterations for convergence and computational time for our numerical examples. When solving the diffusion equation in DSA, we employ an algebraic multigrid  (AMG)  solver implemented based on the iFEM package \cite{Chen:2008ifem}. \pzc{Throughout this section, unless specified, zero initial guesses are employed.}

The online performance of the proposed method is tested with a test set of parameters $\mathcal{P}_{\textrm{test}}$ that lies outside the training set $\mathcal{P}_{\textrm{train}}$. We define the average number of transport sweeps for the test set as 
\begin{equation*}
    \bar{n}_{\textrm{sweep}} = \frac{\sum_{\bmu\in\mathcal{P}_{\pzc{\textrm{test}}}}\textrm{number of transport sweeps for  convergence for the parameter  } \bmu }{\textrm{total number of parameters in }\mathcal{P}_{\textrm{test}}}.
\end{equation*}
In high dimensions, when many angular directions are used in the angular discretization, the computational cost of applying DSA or ROMSAD preconditioners is marginal compared to transport sweeps. Thus, the average computational time is approximately proportional to $\bar{n}_{\textrm{sweep}}$.  To verify that iterative solvers converge to the correct solution, we compute the average $l_\infty$ residual: 
 \begin{equation*}
     \bar{\mathcal{R}}_{\infty} = \frac{\sum_{\bmu\in\mathcal{P}_{\textrm{test}}}||(\BI-\BT\BSigma_s)\brho_{\bmu}-{\wbb}_{\bmu}||_\infty}{\textrm{total number of parameters in }\mathcal{P}_{\textrm{test}}},
 \end{equation*}
where $\brho_{\bmu}$ is the numerical solution obtained by the underlying iterative solver. We define the average relative computational time $\bar{T}_{\textrm{rel}}$ with respect to GMRES-DSA.

\begin{rem}\label{rem:gmres}
We want to point out that when a zero initial guess is used in FGMRES (see Alg. \ref{alg:fgmres}), the number of transport sweeps required for convergence is equal to one plus the number of iterations, as one transport sweep is needed to compute the right-hand side of the memory-efficient formulation, $\wbb$, based on \eqref{eq:memory-efficient-rte}. When other initial guesses are employed, an additional transport sweep will be needed to compute the initial residual. 
\end{rem}

%%%%%%%%%%%%%%%%%%%%%%%%%%%%%%%%%%%%%%%%%%%%%%%%%%%%%%%%%%%%
\subsection{Lattice problem\label{sec:lattice}}
\begin{figure}[]
  \begin{center} 
\includegraphics[width=0.45\textwidth]{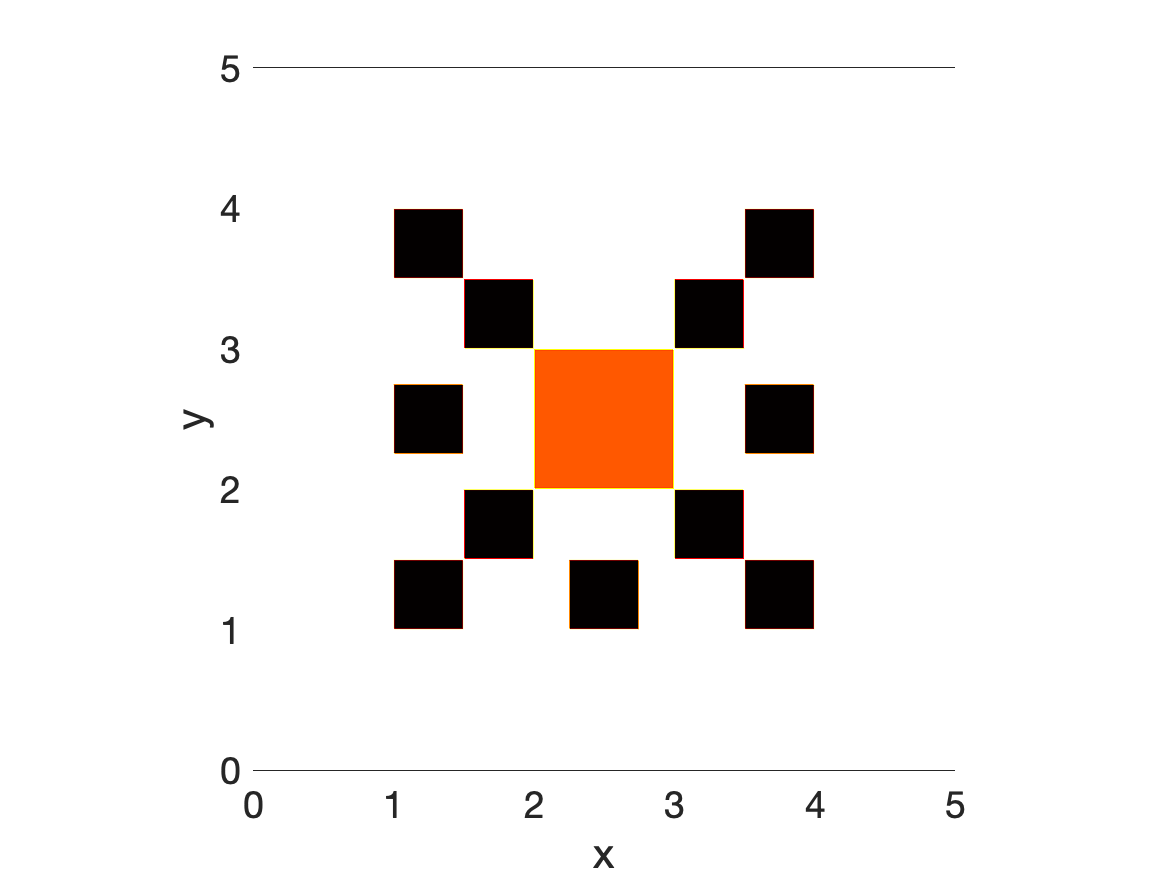}
\includegraphics[width=0.45\textwidth]{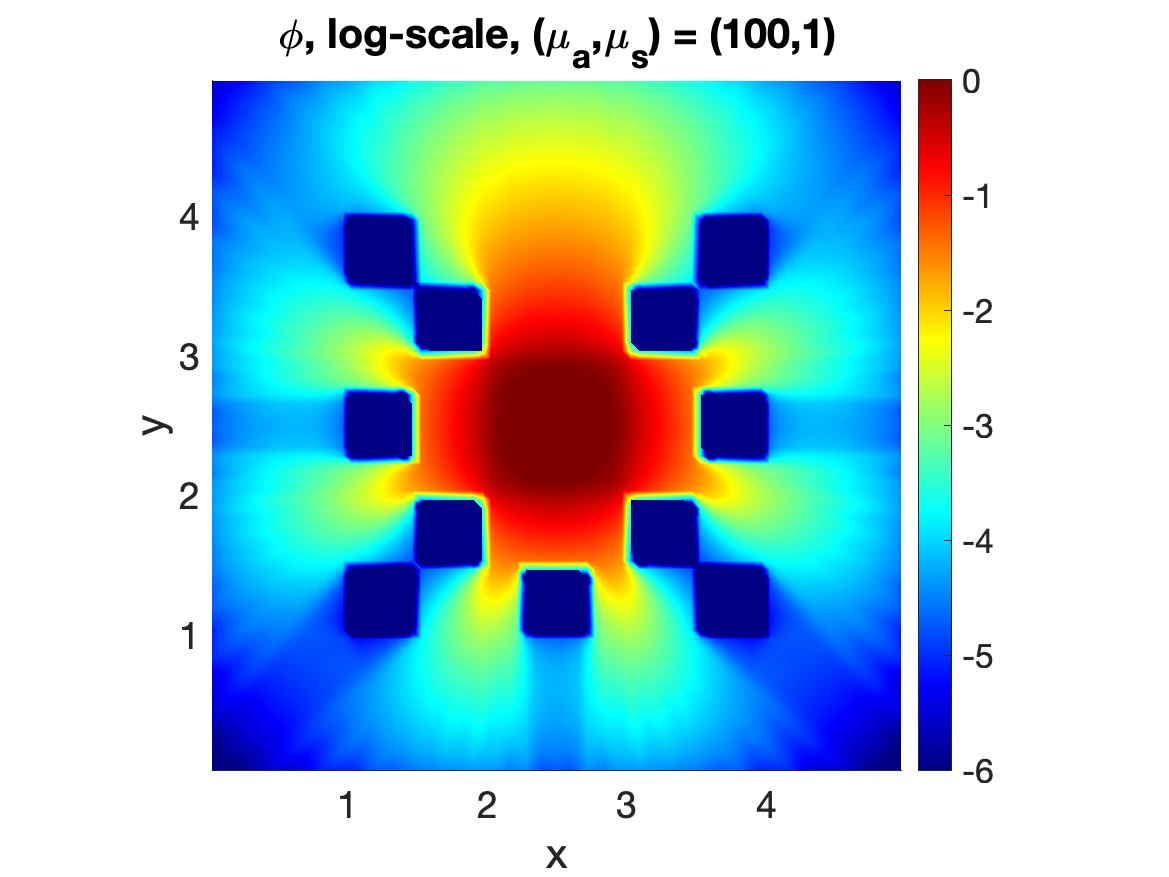}
  \caption{The set-up and a reference solution for the lattice problem in Sec. \ref{sec:lattice}. Left: the set-up for the lattice problem. Black: pure absorption regions with $(\sigma_a,\sigma_s)=(\mu_a,0)$. White and orange: pure scattering regions with $(\sigma_a,\sigma_s)=(0,\mu_s)$. Orange: constant source term with $G=1$. Right: the reference solution under log-scale for $(\mu_a,\mu_s)=(100,1.0)$. We want to point out that negative scalar flux can be generated in this example, but it will not break the linear solver, since we are not considering thermal radiation here. When generating the plot under  log-scale, we take $\max(10^{-16},\phi(\bx))$.  \label{fig:lattice}}
  \end{center}
\end{figure}

We consider a parametric lattice problem with zero inflow boundary conditions in the computational domain $\Gamma_{\bx}=[0,5]^2$. The setup is presented in the left picture of Fig. \ref{fig:lattice}, where the black regions are pure absorption regions  with $(\sigma_a,\sigma_s)=(\mu_a,0)$ and other regions are pure scattering regions with $(\sigma_a,\sigma_s)=(0,\mu_s)$. A source term defined by
\begin{align*}
G(x,y) = \begin{cases}
         1.0,\quad\text{if}\quad |x-2.5|<0.5\;\text{and}\quad|y-2.5|<0.5,\\
         0, \quad\text{otherwise},
         \end{cases}
\end{align*}
is imposed in the orange region at the center of the computational domain. The parameter for this problem $\bmu=(\mu_a,\mu_s)\in[95,105]\times[0.5,1.5]$ determines the strength of absorption  and scattering in the corresponding regions. 
We use an $80\times 80$ uniform mesh in space and CL($40,6$) angular discretization. A reference full order solution is presented in the right picture of Fig. \ref{fig:lattice}. We set the tolerance for the relative residual in (F)GMRES solvers to $10^{-11}$. 

The training set for this problem are $121$ pairs of uniformly sampled $(\mu_a,\mu_s)$ from $[95,105]\times[0.5,1.5]$. We select $\bmu_1=(100,1)$ to start the greedy algorithm.
We test the performance of our method with $10$ pairs of randomly sampled parameters outside the training set. 

\textbf{Online efficiency:} In Tab. \ref{tab:lattice-online}, we present the results using ROMs constructed with $\mfw=2$ and various values of the tolerance $\epsrom$. Compared to GMRES-DSA, FGMRES-ROMSAD achieves comparable accuracy and approximately $2.43$ times acceleration with $\epsrom=10^{-7}$, $2.61$ times acceleration with $\epsrom=10^{-8}$ and $3.00$ times acceleration with $\epsrom=10^{-9}$. As expected, the smaller tolerance $\epsrom$, the greater acceleration is achieved. The average relative computational time is 
approximately proportional to the average number of transport sweeps needed for convergence. 

\textbf{\pzc{Comparison with interpolated initial guesses.}} \pzc{In the offline stage of ROM construction, the greedy algorithm generates high-fidelity solutions for scalar fluxes at sampled parameter pairs. In the online stage, we can also interpolate these high fidelity solutions to generate improved initial guesses for GMRES-DSA. Here, due to scattered nature of these sampled parameters, we employ radial basis function (RBF) interpolation. However, as shown in Tab. \ref{tab:lattice-online}, such interpolated initial guesses provide only marginal acceleration for GMRES-DSA, while FGMRES-ROMSAD achieves significant acceleration. The key advantage of projection-based ROMs over direct interpolation lies in their ability to achieve higher accuracy by explicitly exploiting the information from the original PDE through the projection.}

%%%%%%%%%%%%%%%%%%%%%%%%%%%%%%%%%%%%%%%%%%%%%%%%%%%%%%%%%%%%
%%%%%%%%%%%%%%%%%%%%%%%%%%%%%%%%%%%%%%%%%%%%%%%%%%%%%%%%%%%%
\begin{table}[htbp]
\centering
   \begin{tabular}{|l|c|c|c|c|c|c|c|c|c|c|c|}
    \hline
%%%%%%%%%%%%%%%%%%%%%%%%%%%%%%%%%%%%%%%%%%%%    
$\epsrom$ & & GMRES-DSA  & \pzc{GMRES-DSA-Interp-IG} & FGMRES-ROMSAD\\ \hline
\multirow{3}{*}{$10^{-7}$} & $\bar{n}_{\textrm{sweep}}$ &   $10.7$ & \pzc{$9.4$} & $5.0$  \\ \cline{2-5}
 & $\bar{T}_{\textrm{rel}}$  & $100\%$  & \pzc{$86.30\%$} &  $43.02\%$ \\ \cline{2-5}
 & $\bar{\mathcal{R}}_{\infty}$   & $1.01\times10^{-12}$  & \pzc{$1.20\times10^{-12}$} & $1.49\times 10^{-12}$ \\ \hline
%%%%%%%%%%%%%%%%%%%%%%%%%%%%%%%%%%%%%%%%%%%%
 \multirow{3}{*}{$10^{-8}$} & $\bar{n}_{\textrm{sweep}}$ &   $10.7$ &\pzc{$9.4$} & $4.5$    \\ \cline{2-5}
 & $\bar{T}_{\textrm{rel}}$   & $100\%$ & \pzc{$86.30\%$} &  $38.30\%$  \\ \cline{2-5}
 & $\bar{\mathcal{R}}_{\infty}$  & $1.01\times10^{-12}$ & \pzc{$1.50\times10^{-12}$} & $1.50\times 10^{-12}$  \\ \hline
%%%%%%%%%%%%%%%%%%%%%%%%%%%%%%%%%%%%%%%%%%%%
\multirow{3}{*}{$10^{-9}$} &  $\bar{n}_{\textrm{sweep}}$ &   $10.7$ & \pzc{$9.3$} & $3.9$    \\  \cline{2-5}
 & $\bar{T}_{\textrm{rel}}$  & $100\%$ & \pzc{$86.70\%$} &  $33.30\%$  \\ \cline{2-5}
 & $\bar{\mathcal{R}}_{\infty}$ & $1.01\times10^{-12}$ & \pzc{$1.77\times10^{-12}$}  & $1.55\times 10^{-12}$ \\ \hline
%%%%%%%%%%%%%%%%%%%%%%%%%%%%%%%%%%%%%%%%%%%%
 \end{tabular}
     \caption{Results for the lattice problem in Sec. \ref{sec:lattice} with the window size $\mfw=2$. \pzc{GMRES-DSA-Interp-IG: GMRES-DSA using interpolated initial guesses based on high-fidelity solutions for sampled parameters.} Dimensions of the reduced order spaces: $r=34$ for $\epsrom=10^{-7}$; $r=44$ for $\epsrom=10^{-8}$; $r=50$ for $\epsrom=10^{-9}$.\label{tab:lattice-online}}
\end{table}

\textbf{Offline efficiency:} In Tab. \ref{tab:lattice-offline}, we present the offline computational times for various steps in the greedy algorithm, and compare the computational time of the greedy algorithm $T_{\textrm{greedy}}$ and generating snapshots for all the $121$ parameters in the training set $T_{\textrm{All Snap.}}$. The ratio $T_{\textrm{All Snap.}}/T_{\textrm{greedy}}$ is approximately $4.06$ with $\epsrom=10^{-7}$, $3.37$ with $\epsrom=10^{-8}$ and $3.04$ with $\epsrom=10^{-9}$, demonstrating the computational saving gained by applying the greedy algorithm. The lower the value of $\epsrom$, the more parameters are sampled, leading to the less saving gained offline, but the better efficiency online. Regardless of the  value for $\epsrom$, the main computational cost in the greedy algorithm is consistently from the step 1 and the step 2(a), which are generating right-hand sides of the memory-efficient formulations for all training parameters and high-fidelity linear solves for sampled parameters, respectively.

%%%%%%%%%%%%%%%%%%%%%%%%%%%%%%%%%%%%%%%%%%%%%%%%%%%%%%%%%%%%
\begin{table}[htbp]
    \centering
     %%%%%%%%%%%%%%%%%%%%%%
    \begin{subtable}{\textwidth}
    \centering
    \begin{tabular}{|l|c|c|c|c|c|c|c|c|c|c|c|}
    \hline 
    $T_{\textrm{Step 1}}/T_{\textrm{greedy}}$ & $T_{\textrm{Step 2(a)}}/T_{\textrm{greedy}}$ & $T_{\textrm{Step 2(b)}}/T_{\textrm{greedy}}$ & $T_{\textrm{Step 2(c)}}/T_{\textrm{greedy}}$ & $T_{\textrm{Step 3}}/T_{\textrm{greedy}}$ & $T_{\textrm{All Snap.}}/T_{\textrm{greedy}}$\\ \hline
    $30.55\%$ &  $51.46\%$ & $8.84\%$ & $0.98\%$ & $8.17\%$ & $4.06$
    \\ 
    \hline
    \end{tabular}
    \caption{\normalsize$\epsrom=10^{-7}$, $17$ sampled parameters}
    \end{subtable} 
    \smallskip
    \vspace{-2.5mm}
    %\smallskip
    %%%%%%%%%%%%%%%%%%%%%%%%%%%
    %%%%%%%%%%%%%%%%%%%%%%
    \begin{subtable}{\textwidth}
    \centering
    \begin{tabular}{|l|c|c|c|c|c|c|c|c|c|c|c|}
    \hline 
    $T_{\textrm{Step 1}}/T_{\textrm{greedy}}$ & $T_{\textrm{Step 2(a)}}/T_{\textrm{greedy}}$ & $T_{\textrm{Step 2(b)}}/T_{\textrm{greedy}}$ & $T_{\textrm{Step 2(c)}}/T_{\textrm{greedy}}$ & $T_{\textrm{Step 3}}/T_{\textrm{greedy}}$ & $T_{\textrm{All Snap.}}/T_{\textrm{greedy}}$\\ \hline
    $25.31\%$ &  $55.09\%$ & $9.50\%$ & $1.35\%$ & $8.75\%$ & $3.37$
    \\ 
    \hline
    \end{tabular}
    \caption{\normalsize $\epsrom=10^{-8}$, $22$ sampled parameters}
    \end{subtable}
    \smallskip
    \vspace{-2.5mm}
    %%%%%%%%%%%%%%%%%%%%%%%%%%%
    %%%%%%%%%%%%%%%%%%%%%%
    \begin{subtable}{\textwidth}
    \centering
    \begin{tabular}{|l|c|c|c|c|c|c|c|c|c|c|c|}
    \hline 
    $T_{\textrm{Step 1}}/T_{\textrm{greedy}}$ & $T_{\textrm{Step 2(a)}}/T_{\textrm{greedy}}$ & $T_{\textrm{Step 2(b)}}/T_{\textrm{greedy}}$ & $T_{\textrm{Step 2(c)}}/T_{\textrm{greedy}}$ & $T_{\textrm{Step 3}}/T_{\textrm{greedy}}$ & $T_{\textrm{All Snap.}}/T_{\textrm{greedy}}$\\ \hline
    $22.89\%$ &  $56.79\%$ & $9.77\%$ & $1.57\%$ & $8.98\%$ & $3.04$
    \\ 
    \hline
    \end{tabular}
\caption{\normalsize$\epsrom=10^{-9}$, $25$ sampled parameters}
    \end{subtable}
    %%%%%%%%%%%%%%%%%%%%%%%%%%%
    \caption{Offline results for the lattice problem in Sec. \ref{sec:lattice}. $T_{\textrm{greedy}}$: the computational time of the greedy algorithm offline. $T_{\textrm{Step $k$}}$: the computational time of the $k$-th step in the greedy algorithm offline. $T_{\textrm{All Snap.}}$: the computational time of generating all the snapshots for the $121$ parameters in the training set $\mathcal{P}_{\textrm{train}}$.\label{tab:lattice-offline} }
\end{table}
%%%%%%%%%%%%%%%%%%%%%%%%%%%%%%%%%%%

%%%%%%%%%%%%%%%%%%%%%%%%%%%%%%%%%%%
\subsection{Pin-cell problem \label{sec:pin-cell}}
\begin{figure}[]
  \begin{center} 
\includegraphics[width=0.45\textwidth]{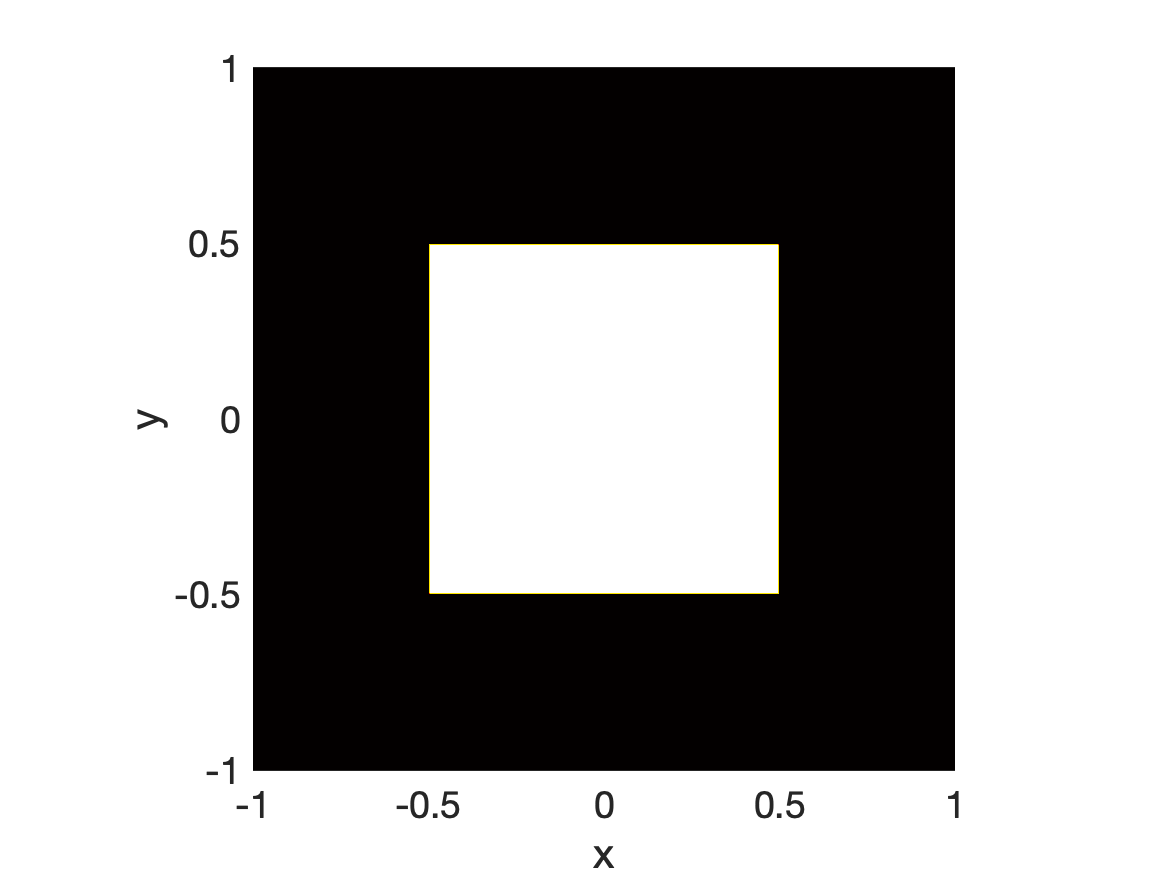}
\includegraphics[width=0.45\textwidth]{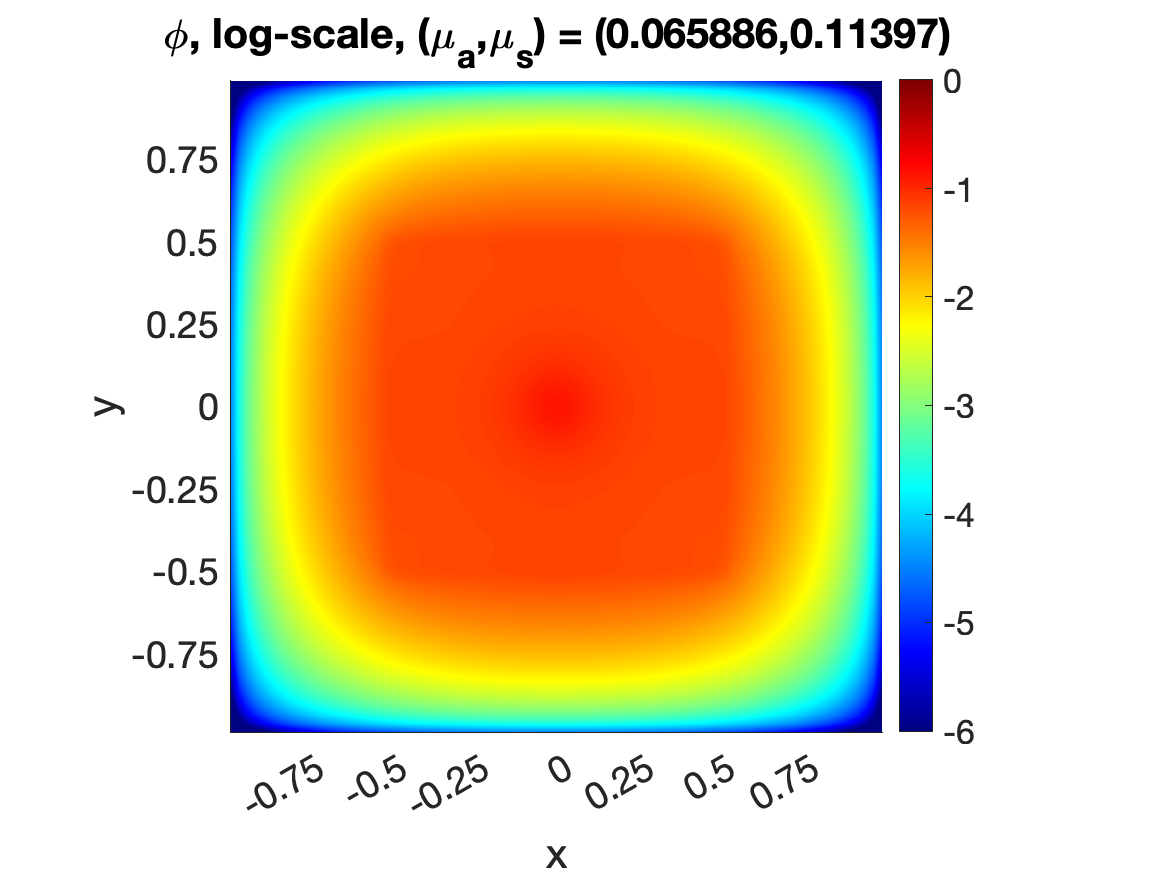}
  \caption{The set-up and reference solutions for the pin-cell problem in Sec. \ref{sec:pin-cell}. Left: problem set-up. Right: $\phi$ for  $(\mu_s,\mu_a)=(0.065886,0.11397)$ (under log-scale). \label{fig:pin-cell}}
  \end{center}
\end{figure}

We consider a parametric pin-cell problem with zero inflow boundary conditions and Gaussian source $G(x,y) = \exp(-100(x^2+y^2))$ on the computational domain $\Gamma_{\bx}=[-1,1]^2$. The  setup is illustrated in the left picture of Fig. \ref{fig:pin-cell}. The outer black region is defined as $\{(x,y):\;|x|>0.5\;\text{or}\;|y|>0.5\}$, where only strong scattering effect ($\sigma_s=100$) exists. The inner white region, defined as $\{(x,y):\;|x|\leq0.5\;\text{and}\;|y|\leq0.5\}$, has parametric scattering and absorption cross sections  given by $(\sigma_a,\sigma_s)=(\mu_a,\mu_s)\in[0.05,0.5]^2$.  The scattering effect in the outer black region is $200$ to $2000$ times as strong as the scattering effect in the inner white region, leading to significant multiscale effects in this problem. We partition the computational domain with an $80\times80$ uniform mesh and apply CL($30,6$) angular discretization. The tolerance for the relative residual of (F)GMRES is set to $10^{-11}$.

The training set for this problem consists of $25$ uniformly sampled parameters from $[0.05,0.5]^2$. The ROM is generated with the initial sample $\bmu_1=(0.275,0.275)$, a window size $\mfw=2$ and a tolerance of $\epsrom=10^{-10}$. We test the performance of the ROMSAD preconditioner with $10$ randomly sampled parameters outside the training set.

\textbf{Online efficiency:}  As shown in Tab. \ref{tab:pin-cell-online}, FGMRES-ROMSAD, on average, converges with $5.1$ transport sweeps per parameter, and achieves approximately $3.33$ times acceleration over GMRES-DSA 
\pzc{using the same zero initial guesses and $3.19$ times acceleration over GMRES-DSA using RBF interpolated initial guesses based on high fidelity solutions for sampled parameters}.

\textbf{Comparison with SI-ROMSAD:}
As presented in \cite{peng2024romsad},  for this pin-cell problem, SI-ROMSAD,whose underlying ROM is constructed using the POD method with data from all the $25$ training parameters and a singular value decomposition (SVD) with truncation tolerance $10^{-9}$, achieves $10$ times acceleration over SI-DSA and $4$ times acceleration over GMRES with a left DSA preconditioner. Using such an accurate ROM, SI-ROMSAD is highly effective.

However, we will use the following setup to demonstrate that when the ROM is not highly accurate, FGMRES-ROMSAD is more robust than SI-ROMSAD. Since the corrections in the first iteration are the same for FGMRES and SI, up to a scaling factor, we set the window size  $\mfw=1$ and $\epsrom=10^{-7}$ when building the ROM with the greedy algorithm. 
This ROM is employed in both FGMRES-ROMSAD and SI-ROMSAD, and we test their performance with the most challenging parameter for ROMSAD in our test set $(\mu_a,\mu_s)=(0.065886,0.11397).$ The full order reference solution for this parameter is presented in the right picture of Fig. \ref{fig:pin-cell}. Sharp features in the scalar flux can be observed near the material interface. The stopping criterion for SI is  whether the absolute $l_2$ residual $||\br^{(l)}||_2=||\bphi^{(l,*)}-\bphi^{(l-1)}||_2$  is smaller than $10^{-11}\times||\wbb||_2$ with $\wbb$ being the right hand side for the memory-efficient formulation. 

\begin{figure}[]
  \begin{center} 
\includegraphics[width=0.45\textwidth]{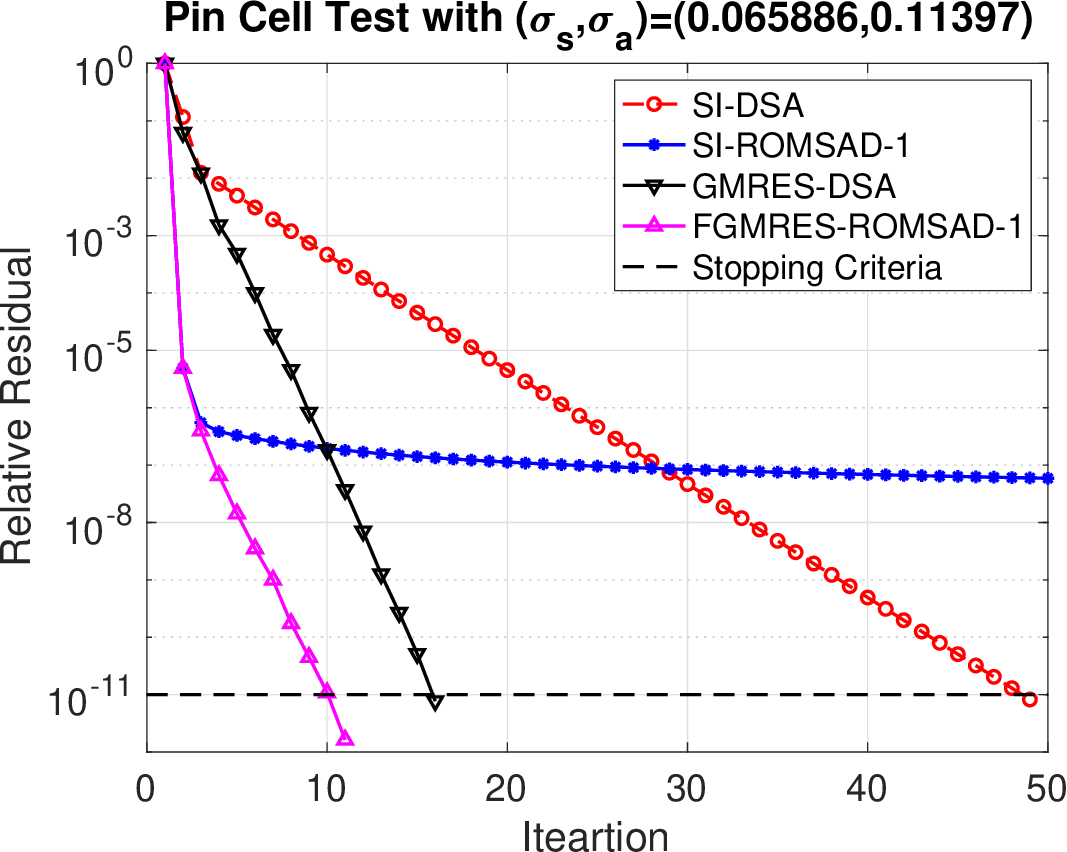}
  \caption{\pzc{Convergence history for FGMRES-ROMSAD and SI-ROMSAD with $\mfw=1$ and $\epsrom=10^{-7}$ for the pin-cell problem in Sec. \ref{sec:pin-cell} with $(\mu_s,\mu_a)=(0.065886,0.11397)$ (under log-scale).} \label{fig:pin-cell-gmres-comparison}}
  \end{center}
\end{figure}

As shown in Tab. \ref{tab:pin-cell-compare} \pzc{and Fig. \ref{fig:pin-cell-gmres-comparison}}, even with this relatively less accurate ROM, FGMRES-ROMSAD still achieves $1.33$ times acceleration over GMRES-DSA and $4.71$ times acceleration over SI-DSA.  In contrast,  SI-ROMSAD is unable to reach the desired accuracy within $100$ iterations, while SI-DSA converges in $49$ iterations. In both FGMRES-ROMSAD and SI-ROMSAD, DSA starts to be applied after the first SA correction. \pzc{As demonstrated in Fig. \ref{fig:pin-cell-gmres-comparison}}, the relative $l_2$ residuals after the first correction of FGMRES-ROMSAD and SI-ROMSAD are comparable and much smaller than that of GMRES-DSA and SI-DSA, demonstrating the effectiveness of the ROM-based correction. However, despite using the same preconditioner from the second iteration and starting from a smaller relative residual, SI-ROMSAD surprisingly converges much more slowly than SI-DSA. We suspect that this slowdown is due to the robustness degeneration of SI-DSA for problems with discontinuous materials \cite{azmy2002unconditionally,warsa2004krylov}.

\textbf{Offline efficiency:} The training set consists of a total of $25$ parameters. In the offline stage, the greedy algorithm samples 22 parameters from this set.  As shown in Tab. \ref{tab:pin-cell-offline}, the computational time of the greedy algorithm is nearly the same as generating high-fidelity snapshots for the entire training set. In summary, when the number of representative parameters is close to the size of the training set, the benefit in using the greedy offline algorithm is less pronounced. 
%%%%%%%%%%%%%%%%%%%%%%%%%%%%%%%%%%%%%%%%%%%%%

% %%%%%%%%%%%%%%%%%%%%%%%%%%%%%%%%%%%%%%%%%%%%%%%%%%%%%%%%%%%%
\begin{table}[htbp]
\centering
   \begin{tabular}{|l|c|c|c|c|c|c|c|c|c|c|c|}
    \hline
%%%%%%%%%%%%%%%%%%%%%%%%%%%%%%%%%%%%%%%%%%%%    
&   GMRES-DSA & \pzc{GMRES-DSA-Interp-IG}  & FGMRES-ROMSAD \\ \hline
 $\bar{n}_{\textrm{sweep}}$ &   $15.5$& \pzc{14.2} & $5.1$  \\ \hline
 $\bar{T}_{\textrm{rel}}$ &   $100\%$ & \pzc{$95.75\%$} & $30.05\%$ \\ \hline
 $\bar{\mathcal{R}}_{\infty}$ &   $5.14\times10^{-13}$ & \pzc{$5.26\times10^{-13}$} & $5.75\times10^{-13}$  \\ \hline
%%%%%%%%%%%%%%%%%%%%%%%%%%%%%%%%%%%%%%%%%%%%
 \end{tabular}
     \caption{Online results for the pin-cell problem in Sec. \ref{sec:pin-cell}. \pzc{GMRES-DSA-Interp-IG: GMRES-DSA using interpolated initial guesses based on high-fidelity solutions for sampled parameters.} In the offline stage, $22$ parameters are sampled, resulting in a  reduced order space whose dimension is $44$. 
     \label{tab:pin-cell-online}}
\end{table}

%%%%%%%%%%%%%%%%%%%%%%%%%%%%%%%%%%%%%%%%%%%%%

%%%%%%%%%%%%%%%%%%%%%%%%%%%%%%%%%%%%%%%%%%%%%%%%%%%%%%%%%%%%
\begin{table}[htbp]
\centering
   \begin{tabular}{|l|c|c|c|c|c|c|c|c|c|c|c|}
    \hline
%%%%%%%%%%%%%%%%%%%%%%%%%%%%%%%%%%%%%%%%%%%%    
&   GMRES-DSA & SI-DSA  & FGMRES-ROMSAD & SI-ROMSAD \\ \hline
 $n_{\textrm{sweep}}$ &   $15$& $49$ & $10$ & $100$ ($\star$)  \\ \hline
 $T_{\textrm{rel}}$ &   $100\%$ & $323.05\%$& $68.55\%$ & $663.47\%$\\ \hline
 $\mathcal{R}_{\infty}$ &  $7.72\times10^{-13}$ & $4.35\times10^{-13}$ &$2.20\times10^{-13}$ & $1.18\times10^{-9}$  \\ \hline
$\mathcal{R}_{\textrm{rel.},2}$ after the first correction & $6.23\times10^{-2}$ & $1.17\times10^{-1}$ & $4.87\times10^{-6}$ & $4.90\times10^{-6}$ 
 \\ \hline
%%%%%%%%%%%%%%%%%%%%%%%%%%%%%%%%%%%%%%%%%%%%
 \end{tabular}
     \caption{Comparision between FGMRES-ROMSAD and SI-ROMSAD for the pin-cell problem in Sec. \ref{sec:pin-cell} with $(\mu_a,\mu_s)=(0.06589,0.11397)$. In the offline stage, $\epsrom=10^{-7}$, $\mfw=1$ and $13$ parameters are sampled, resulting in a  reduced order space whose dimension is $13$. $\mathcal{R}_{\textrm{rel.},2}$: relative $l_2$ residual.
     $\star$: means that SI-ROMSAD fails to converge as its relative $l_2$ residual almost gets stuck and decays very slowly. 
     \label{tab:pin-cell-compare}}
\end{table}
%%%%%%%%%%%%%%%%%%%%%%%%%%%%%%%%%%%%%%%%%%%%%%%%%%%%%%%%%%%%

%%%%%%%%%%%%%%%%%%%%%%
\begin{table}[htbp]
    \centering
    \begin{tabular}{|l|c|c|c|c|c|c|c|c|c|c|c|}
    \hline 
    $T_{\textrm{Step 1}}/T_{\textrm{greedy}}$ & $T_{\textrm{Step 2(a)}}/T_{\textrm{greedy}}$ & $T_{\textrm{Step 2(b)}}/T_{\textrm{greedy}}$ & $T_{\textrm{Step 2(c)}}/T_{\textrm{greedy}}$ & $T_{\textrm{Step 3}}/T_{\textrm{greedy}}$ & $T_{\textrm{All Snap.}}/T_{\textrm{greedy}}$\\ \hline
    $5.08\%$ &  $83.07\%$ & $9.40\%$ & $1.39\%$ & $1.04\%$ & $1.02$
    \\ 
    \hline
    \end{tabular}
    \caption{Offline results for the pin-cell problem in Sec. \ref{sec:pin-cell}. $T_{\textrm{greedy}}$: the computational time of the greedy algorithm offline. $T_{\textrm{Step $k$}}$: the computational time of the $k$-th step in the greedy algorithm offline. $T_{\textrm{All Snap.}}$: the computational time of generating all the snapshots for the $51$ parameters in the training set $\mathcal{P}_{\textrm{train}}$ using GMRES-DSA.\label{tab:pin-cell-offline} }
\end{table}

%%%%%%%%%%%%%%%%%%%%%%%%%%%%%%%%%%%%%%%%%%%%%%%%%%%%%%%%%%%%
\subsection{Variable scattering problem\label{sec:variable-scattering}}
\begin{figure}[]
  \begin{center} 
\includegraphics[width=0.45\textwidth]{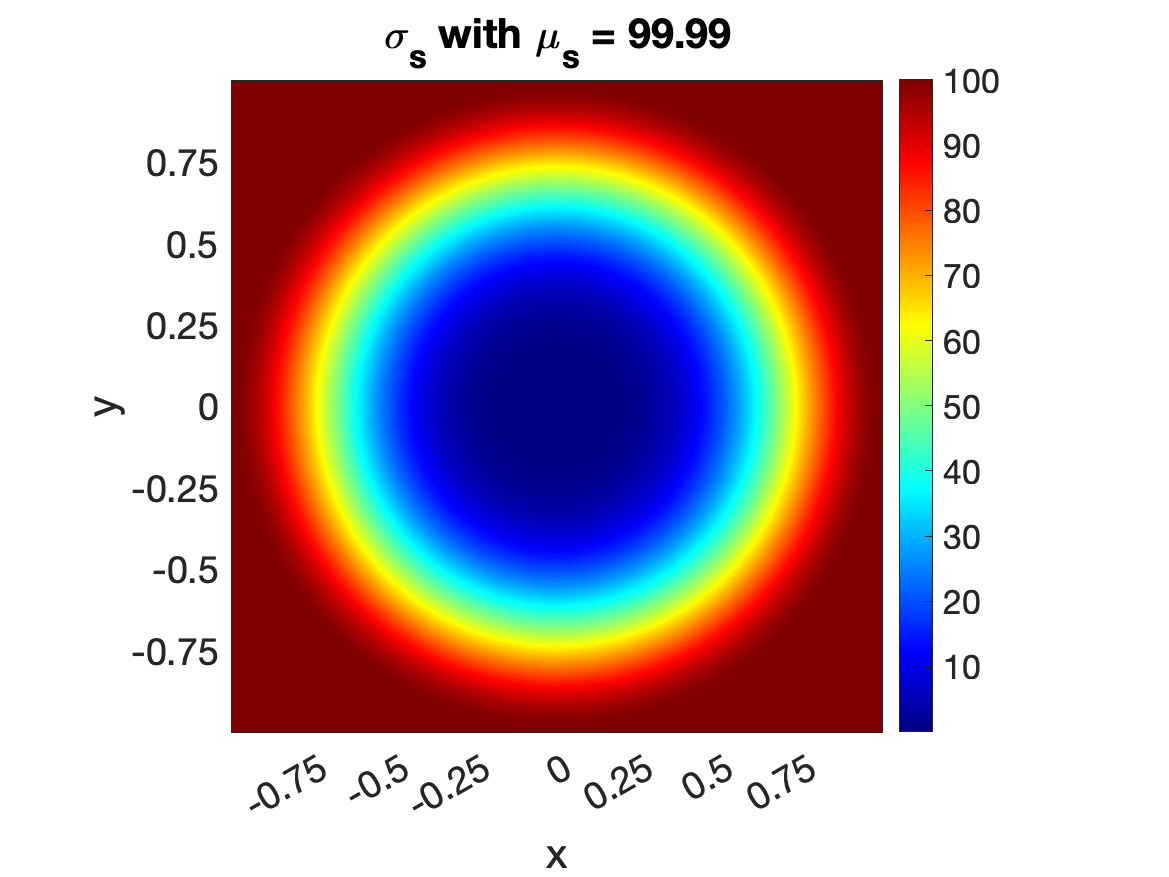}
\includegraphics[width=0.45\textwidth]{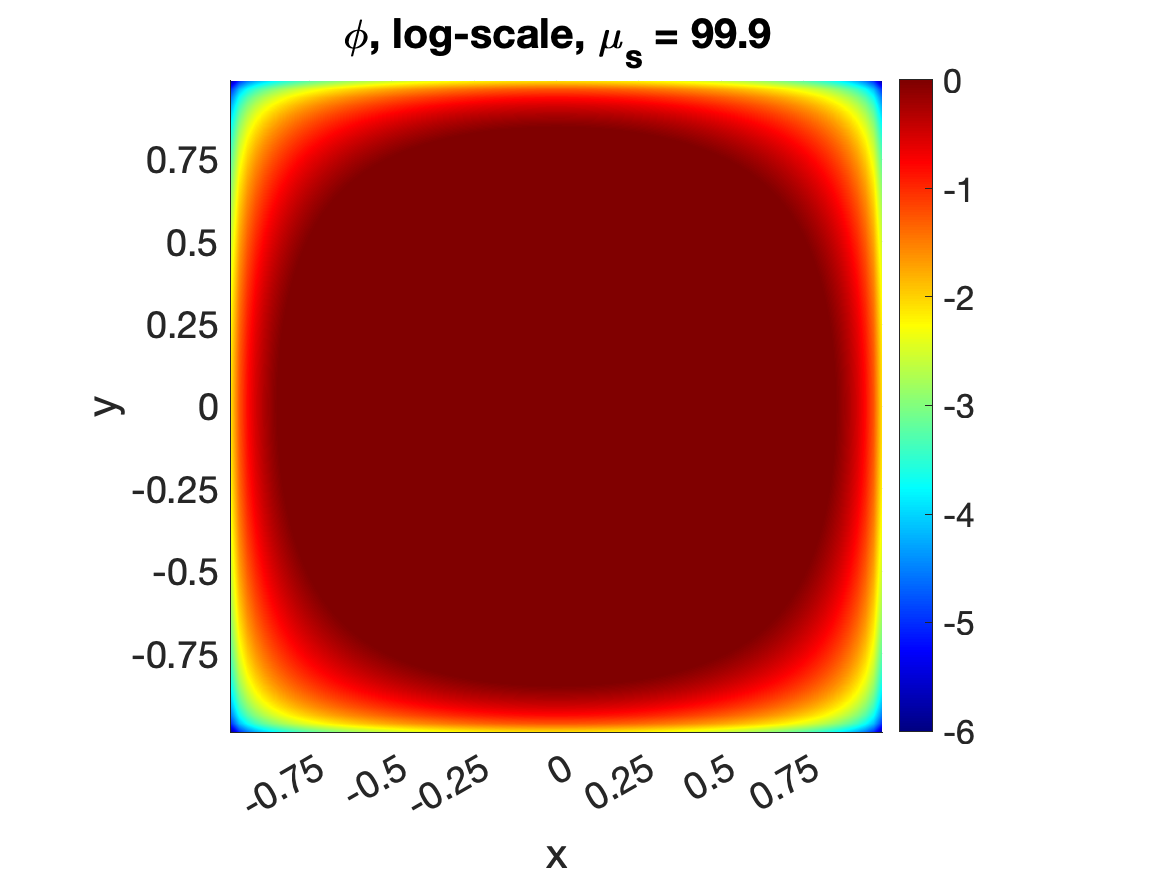}
  \caption{Left: $\sigma_s(\bx)$ with $\mu_s=99.9$ for the variable scattering problem in Sec. \ref{sec:variable-scattering}. Right: corresponding angular flux. \label{fig:variable-scattering}}
  \end{center}
\end{figure}

The final test involves a parametric variable scattering problem on the computational domain  $\Gamma_{\bx}=[-1,1]^2$. The parametric scattering cross section is defined as
\begin{equation}
    \sigma_s(x,y) = \begin{cases}
        \mu_s r^4(2-r^2)^2+0.1,\quad r=\sqrt{x^2+y^2}\leq 1,\\
        \mu_s+0.1,\quad\text{otherwise}, 
    \end{cases}\quad\mu_s\in[49.9,99.9],\quad\sigma_a(x,y)=0.
\end{equation}
 The scattering effect smoothly changes from $0.1$ to $\mu_s+0.1\in [50,100]$ when moving from the center of the computational domain to its boundary.  In other words, there is a smooth transition from transport dominant regime to scattering dominant regime in this problem. The parameter $\mu_s$ determines the speed of this transition.
The boundary condition is zero inflow boundary conditions. A Gaussian source 
 $G(x,y)=\frac{10}{\pi}\exp(-100(x^2+y^2))$ is imposed within the computational domain. 
 An  $80\times 80$ uniform spatial mesh  and CL($30,6$) angular discretization are applied. In Fig. \ref{fig:variable-scattering}, we present 
a configuration of the scattering cross section $\sigma_s(\bx)$ and the corresponding high-fidelity reference solution. 

The parameter of this problem, $\mu_s\in[49.9,99.9]$, determines how quickly the scattering strength increases from the center to the boundary of the computational domain. We use a training set consisting of $51$ uniformly sampled parameters, and set the initial sample, the window size and the tolerance in the greedy algorithm as $\mu_{s,1}=74.9$, $\mfw=2$ and $\epsrom=10^{-9}$. The tolerance of relative residual for (F)GMRES solver is set to $10^{-11}$. We use $10$ randomly sampled parameters outside the training set to test the performance of our ROMSAD preconditioner.

%%%%%%%%%%%%%%%%%%%%%%%%%%%%%%%%%%%%%%%%%%%%%%

%%%%%%%%%%%%%%%%%%%%%%%%%%%%%%%%%%%%%%%%%%%%%%%%%%%%%%%%%%%%
\textbf{Online efficiency:}  In this example, we test the performance of ROMSAD coupled with different DSA preconditioners. Both the fully consistent (FC) and partially consistent (PC) DSA preconditioners are considered (see  \ref{sec:dsa}
 for details).

As shown in Tab. \ref{tab:variable-scattering-online}, in GMRES-DSA, the right FC-DSA preconditioner results in a $26.67\%$ reduction in the number of transport sweeps and $9.35\%$ reduction in computational time compared to the right PC-DSA preconditioner. We suspect that the different percentage of the reduction in the number of transport sweeps and the computational time is due to the more computational cost required by FC-DSA in comparison to PC-DSA. Conversely, by using ROM-based corrections in the first two iterations, FGMRES-ROMSAD converges in an average of $4.7$ transport sweeps when using PC-DSA  and $4.6$ transport sweeps when using FC-DSA. Coupling ROM-based corrections with both FC-DSA and PC-DSA, FGMRES achieves at least $3.11$ times acceleration over GMRES with PC-DSA preconditioner and $2.82$ times acceleration over GMRES with FC-DSA preconditioner. 

\pzc{Using interpolated initial guesses, GMRES-PC-DSA and GMRES-FC-DSA converge in an average of $7.8$ and $6.3$ transport sweeps, respectively. This results in approximately $1.92$ and $1.79$ times acceleration for GMRES-PC-DSA and GMRES-FC-DSA with zero initial guesses, respectively.  However, this acceleration is still less significant than proposed FGMRES-ROMSAD method which results in $3.1$ and $2.82$ times acceleration for PC-DSA and FC-DSA. However, this acceleration remains less significant than that achieved by the proposed FGMRES-ROMSAD method. For brevity, we omit the details of GMRES-DSA with interpolated initial guesses in Tab. \ref{tab:variable-scattering-online}.}

\textbf{Offline efficiency:} We observe that whether generating snapshots with the partially or fully consistent DSA preconditioner, the same parameters are sampled, and the difference between the resulting reduced order basis is negligible. Therefore, we present only the offline results with the fully consistent DSA in Tab. \ref{tab:variable-scattering-offline}. The main cost of the greedy algorithm  still comes from its step 1 and step 2(a), which are generating the right-hand sides for the memory-efficient formulation and the high fidelity solves for sampled parameters, respectively. Using the greedy algorithm to build the ROM only requires approximately $25.00\%$ of the time needed to generate snapshots for all training parameters using GMRES with right FC-DSA preconditioner.

%%%%%%%%%%%%%%%%%%%%%%%%%%%%%%%%%%%%%%%%%%%%%%%%%%%%%%%%%%%%
\begin{table}[htbp]
\centering
   \begin{tabular}{|l|c|c|c|c|c|c|c|c|c|c|c|}
    \hline
%%%%%%%%%%%%%%%%%%%%%%%%%%%%%%%%%%%%%%%%%%%%    
&   GMRES-PC-DSA & GMRES-FC-DSA & FGMRES-PC-ROMSAD & FGMRES-FC-ROMSAD\\ \hline
 $\bar{n}_{\textrm{sweep}}$ &   $15.0$& $11.0$ & $4.7$ & $4.6$  \\ \hline
 $\bar{T}_{\textrm{rel}}$ &   $100\%$ &  $90.65\%$ & $30.61\%$ & $32.11\%$ \\ \hline
 $\bar{\mathcal{R}}_{\infty}$ &   $9.26\times10^{-13}$ & $3.22\times10^{-13}$ &    $1.05\times10^{-12}$ & $5.32\times10^{-13}$  \\ \hline
%%%%%%%%%%%%%%%%%%%%%%%%%%%%%%%%%%%%%%%%%%%%
 \end{tabular}
     \caption{Online results for the variable scattering problem in Sec. \ref{sec:variable-scattering}. In the offline stage, $8$ parameters are sampled, resulting in a  reduced order space whose dimension is $16$. 
     GMRES-PC/FC-DSA: GMRES with partially/fully consistent right DSA preconditioner. FGMRES-PC/FC-ROMSAD: FGMRES with ROMSAD preconditioner using partially/fully consistent DSA preconditioner.
     The relative computational cost is computed with respect to GMRES-PC-DSA.\label{tab:variable-scattering-online}}
\end{table}

%%%%%%%%%%%%%%%%%%%%%%%%%%%%%%%%%%%%%%%%%%%%%%%%%%%%%%%%%%%%
\begin{table}[htbp]
    \centering
    \begin{tabular}{|l|c|c|c|c|c|c|c|c|c|c|c|}
    \hline 
    $T_{\textrm{Step 1}}/T_{\textrm{greedy}}$ & $T_{\textrm{Step 2(a)}}/T_{\textrm{greedy}}$ & $T_{\textrm{Step 2(b)}}/T_{\textrm{greedy}}$ & $T_{\textrm{Step 2(c)}}/T_{\textrm{greedy}}$ & $T_{\textrm{Step 3}}/T_{\textrm{greedy}}$ & $T_{\textrm{All Snap.}}/T_{\textrm{greedy}}$\\ \hline
    $26.08\%$ &  $61.61\%$ & $8.51\%$ & $0.44\%$ & $3.3\%$ & $4.00$
    \\ 
    \hline
    \end{tabular}
    \caption{Offline results for the variable scattering problem in Sec. \ref{sec:variable-scattering}. $T_{\textrm{greedy}}$: the computational time of the greedy algorithm offline. $T_{\textrm{Step $k$}}$: the computational time of the $k$-th step in the greedy algorithm offline. $T_{\textrm{All Snap.}}$: the computational time of generating all the snapshots for the $51$ parameters in the training set $\mathcal{P}_{\textrm{train}}$ using GMRES with right fully consistent DSA preconditioner.\label{tab:variable-scattering-offline} }
\end{table}

%%%%%%%%%%%%%%%%%%%%%%%%%%%%%%%%%%%%%%%%%%%%%%%%%%%%%%%%%%%%
\subsection{\pzc{A two-material problem with parametric inflow boundary conditions}\label{sec:parametric-bc}}
\begin{figure}[]
  \begin{center} 
\includegraphics[width=0.45\textwidth]{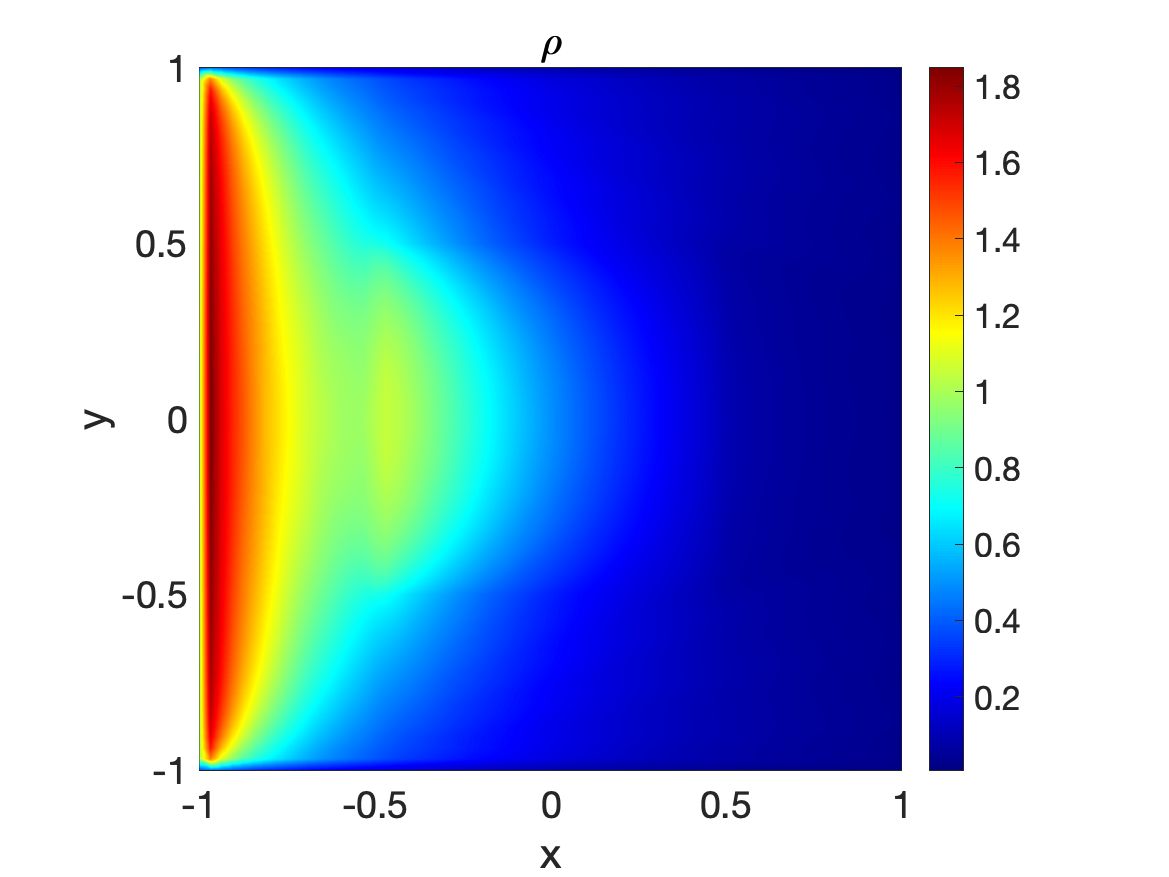}
\includegraphics[width=0.45\textwidth]{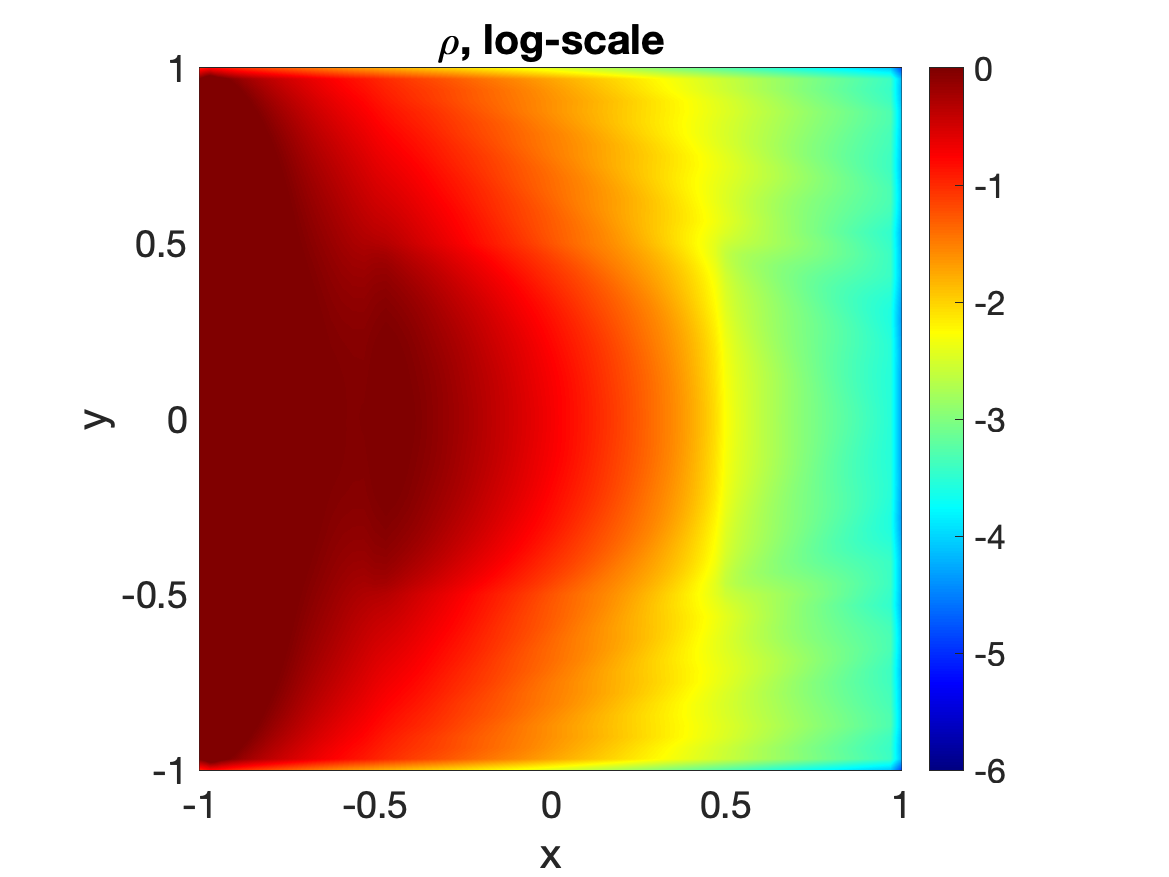}
  \caption{\pzc{Reference solution with $(\mu_s,\mu_{\textrm{bc}})=(3,1)$ for the two-material problem with parametric inflow boundary conditions in Sec. \ref{sec:parametric-bc}.} \label{fig:parametric-bc}}
  \end{center}
\end{figure}

\pzc{We consider a problem with parametric inflow boundary conditions and parametric scattering cross section on the computational domain $\Gamma_{\bx}=[-1,1]^2$. The source is $G(x,y)=0$. The scattering and absorption sections are defined as}
\begin{equation}
\pzc{
\sigma_s(x,y)=\begin{cases}
              \mu_s,\quad\text{where }|x|\leq 0.5\text{ and }|y|\leq 0.5,\\
              0.0,\quad\text{otherwise,}
              \end{cases}
\quad\text{and}\quad
\sigma_s(x,y)=\begin{cases}
              0.0,\quad\text{where }|x|\leq 0.5\text{ and }|y|\leq 0.5,\\
              1.0.\quad\text{otherwise.}
              \end{cases}
}
\end{equation}
\pzc{Zero inflow boundary conditions are imposed to the right, top and bottom boundaries of the computational domain, while the left boundary has a parametric inflow boundary condition 
\begin{equation}
\psi(-1,y,\BOmega_x,\BOmega_y)=\mu_{\textrm{bc}},\quad \text{if}\quad\BOmega_x\geq 0.
\end{equation} 
The parameter of this problem $\bmu=(\mu_s,\mu_{\textrm{bc}})\in[1,5]\times[0.5,1.5]$ determines the scattering strength of the inner cell and the inflow boundary condition on the left boundary. We partition the computational domain with a $64\times64$ uniform mesh and use CL$(60,6)$ angular discretization. In Fig. \ref{fig:parametric-bc}, we present a high-fidelity reference solution. We set the tolerance for the relative residual of (F)GMRES as $10^{-11}$.}

\pzc{The training set of this problem is chosen as the tensor product of $\{\mu_s=1+\frac{j}{2},\;0\leq j\leq 8,\;j\in\mathbb{Z}\}$ and $\{\mu_{\textrm{bc}}=0.5+0.1k,\;0\leq k\leq 10,\;k\in\mathbb{Z}\}$. We select the initial sample of the greedy algorithm as $\bmu_1=(\mu_{s,1},\mu_{\textrm{bc},1})=(3,1)$. The window size and the tolerance in the greedy algorithm are $\mfw=2$ and $\epsrom=10^{-9}$, respectively. We test the performance of our method with $10$ pairs randomly sampled parameters outside the training set.}

\pzc{\textbf{Online efficiency:}} \pzc{As shown in Tab. \ref{tab:parametric-bc-online}, using zero initial guess, FGMRES-ROMSAD, on average, achieves approximately $2.48$ times acceleration over GMRES-DSA.} 

%%%%%%%%%%%%%%%%%%%%%%%%%%%%%%%%%%%%%%%%%%%%%%
%%%%%%%%%%%%%%%%%%%%%%%%%%%%%%%%%%%%%%%%%%%%%%

%%%%%%%%%%%%%%%%%%%%%%%%%%%%%%%%%%%%%%%%%%%%%%%%%%%%%%%%%%%%
\begin{table}[htbp]
\centering
    \begin{tabular}{|l|c|c|c|c|c|c|c|c|c|c|c|}
     \hline
 %%%%%%%%%%%%%%%%%%%%%%%%%%%%%%%%%%%%%%%%%%%%    
 &   \pzc{GMRES-DSA}  & \pzc{FGMRES-ROMSAD} %&  \pzc{GMRES-DSA-ROMIG}
 \\ \hline
  \pzc{$\bar{n}_{\textrm{sweep}}$} &   \pzc{$10.0$} & \pzc{$4.3$} %& \pzc{$4.3$}   
  \\ \hline
  \pzc{$\bar{T}_{\textrm{rel}}$} &   \pzc{$100\%$} &  \pzc{$40.35\%$} %& \pzc{$41.13\%$} 
  \\ \hline
  \pzc{$\bar{\mathcal{R}}_{\infty}$} &   \pzc{$5.88\times10^{-12}$} & \pzc{$2.30\times10^{-12}$} %& \pzc{$2.63\times10^{-12}$}  
  \\ \hline
% %%%%%%%%%%%%%%%%%%%%%%%%%%%%%%%%%%%%%%%%%%%%
  \end{tabular}
      \caption{\pzc{Online results for the two-material problem with parametric boundary conditions in Sec. \ref{sec:parametric-bc}. In the offline stage, $9$ parameters are sampled, resulting in a  reduced order space whose dimension is $18$.} 
      \label{tab:parametric-bc-online}}
 \end{table}

\pzc{\textbf{Offline efficiency:}} \pzc{As presented in Tab. \ref{tab:parametric-bc-offline}, compared with generating snapshots for the entire training set, building with the ROM greedy algorithm only requires approximately $18.14\%$ computational time.} 

%%%%%%%%%%%%%%%%%%%%%%%%%%%%%%%%%%%%%%%%%%%%%%%%%%%%%%%%%%%%
\begin{table}[htbp]
    \centering
    \begin{tabular}{|l|c|c|c|c|c|c|c|c|c|c|c|}
    \hline 
    \pzc{$T_{\textrm{Step 1}}/T_{\textrm{greedy}}$} & \pzc{$T_{\textrm{Step 2(a)}}/T_{\textrm{greedy}}$} & \pzc{$T_{\textrm{Step 2(b)}}/T_{\textrm{greedy}}$} & \pzc{$T_{\textrm{Step 2(c)}}/T_{\textrm{greedy}}$} & \pzc{$T_{\textrm{Step 3}}/T_{\textrm{greedy}}$} & \pzc{$T_{\textrm{All Snap.}}/T_{\textrm{greedy}}$}\\ \hline
    \pzc{$7.76\%$} &  \pzc{$7.75\%$} & \pzc{$1.46\%$} & \pzc{$0.94\%$} & \pzc{$1.09\%$} & \pzc{$5.51$}
    \\ 
    \hline
    \end{tabular}
    \caption{\pzc{Offline results for the two-material problem with parametric inflow boundary conditions in Sec. \ref{sec:parametric-bc}. $T_{\textrm{greedy}}$: the computational time of the greedy algorithm offline. $T_{\textrm{Step $k$}}$: the computational time of the $k$-th step in the greedy algorithm offline. $T_{\textrm{All Snap.}}$: the computational time of generating all the snapshots for the $99$ parameters in the training set $\mathcal{P}_{\textrm{train}}$ using GMRES with right fully consistent DSA preconditioner.}\label{tab:parametric-bc-offline} }
\end{table}

%%%%%%%%%%%%%%%%%%%%%%%%%%%%%%%%%%%%%%%%%%%%%%%%%%%%%%%%%%%%
\section{Conclusions\label{sec:conclusion}}
In this paper, we extend the ROMSAD preconditioner from the SI framework to the Krylov framework and construct the underlying ROM more efficiently through a greedy algorithm.

From our numerical tests, we make the following observations. 
\begin{enumerate}
    \item  In the online stage, FGMRES-ROMSAD outperforms GMRES with the right DSA preconditioner. 
    \item When the underlying ROM is not highly accurate, FGMRES-ROMSAD is more robust than SI-ROMSAD.
    \item In the offline stage, if the number of representative parameters is much smaller than the size of the training set, the greedy algorithm can construct the ROM significantly more efficiently than methods that require high-fidelity data for the entire training set. 
\end{enumerate}

\pzc{Our current method focuses on linear parametric RTE with affine dependence on the underlying parameters. Extending our approach to nonlinear thermal radiation or non-affine parametric problems will require the use of classical hyper reduction techniques, such as EIM \cite{barrault2004empirical} and DEIM \cite{chaturantabut2010nonlinear}, in the construction of the ROM.} \pzc{Additionally}, in the future,  we aim to integrate the proposed method as  a building block for multi-query applications such as uncertainty quantification, inverse problems and design optimization. 

\pzc{Besides the algorithm development, it is essential to theoretically analyze the relation between the converge rate of the Krylov solver with ROM-enhanced preconditioenrs and the Kolmogrov $n$-width of the solution manifold for the underlying parametric problem, i.e. the best error of approximating this manifold with a $n$-dimensional linear space. To the best of our knowledge, such analysis is not available even for Krylov methods with ROM-based preconditioners solving parametric elliptic problems. Convergence theories for ROMs and analysis tools for inexact and flexible Krylov solvers, e.g. \cite{simoncini2003theory}, may be leveraged to establish such theories.}

%%%%%%%%%%%%%%%%%%%%%%%%%%%%%%%%%%%%%%%%%%%%%%%%%%%%%%%%%%%%
\section*{Acknowledgement}
The author would like to thank Prof. Fengyan Li and  from Rensselaer Polytechnic Institute and her student Kimberly Matsuda for their discussions, and reviewers of his previous work \cite{peng2024romsad} for their interesting questions which motivates this paper. \pzc{Additionally, the author would like to thank the reviewers of this draft for their careful reading and valauble feedbacks.}

\section*{CRediT authorship contribution statement}

{\bf Zhichao Peng:} Writing – original draft, Writing – review \& editing, Visualization, Validation, Software, Methodology, Data curation, Conceptualization.

\section*{Declaration of generative AI and AI-assisted technologies in the writing process}
During the preparation of this work the author(s) used ChatGPT in order to check grammar errors and improve readability. After using this tool/service, the author(s) reviewed and edited the content as needed and take(s) full responsibility for the content of the publication.
%%%%%%%%%%%%%%%%%%%%%%%%%%%%%%%%%%%%%%%%%%%%%%%%%%%%%%%%%%%%

\appendix

\section{Partially and fully consistent DSA \label{sec:dsa}}
Following \cite{adams2001discontinuous}, we briefly outline how to derive a fully/partially consistent discretization for DSA. For simplicity, we consider 1D slab geometry with a angular discretization using a quadrature rule $\{(\pm v_j,\omega_{\pm j})\}$ satsfying $\omega_j=\omega_{-j}$.

Let $\{T_i=[x_{i-\half},x_{i+\half}], \; i=1,\dots, N_x\}$ be a partition of the computational domain and $\{\eta_i(\bx)\}_{i=1}^{N_{\textrm{DOF}}}$ be an orthonormal basis for DG discrete space. We introduce the discrete advection operator with central flux $\BD_C\in\mathbb{R}^{N_{\textrm{DOF}}\times N_{\textrm{DOF}}}$ and the jump operator $\BD_J\in\mathbb{R}^{N_{\textrm{DOF}}\times N_{\textrm{DOF}}}$ as:
\begin{subequations}
\begin{align}
&(\BD_C)_{kl}=-\sum_{i=1}^{N_x}\left(\int_{T_i} \partial_x\eta_k(x) \eta_l(x) dx+\{\eta_l(x_{i+\half})\} [\eta_k(x_{i+\half})]\right),\;
(\BD_J)_{kl}=-\sum_{i=1}^{N_x}[\eta_l(x_{i+\half})][\eta_k(x_{i+\half})],\\
&\{\eta_{l}(x_{i+\frac{1}{2}})\}=\left(\eta_{l}(x^+_{i+\frac{1}{2}})+\eta_{l}(x^-_{i+\frac{1}{2}})\right)/2\quad\text{and}\quad
[\phi_k(x_{i+\half})]=\eta_{l}(x^+_{i+\frac{1}{2}})-\eta_{l}(x^-_{i+\frac{1}{2}})
\end{align}
\end{subequations}
Then, the discrete upwind advection operator is
$ \BD_j = v_j\left(\BD_C-\frac{\textrm{sign}(v_j)}{2}D_J\right).$

Assume $\delta \psi^{(l)}(x,v_j)$ can be represented as a $P_1$ expansion $\delta\psi^{(l)}(x,v_j) = \delta\phi^{(l)}(x)+3v_j \delta J^{(l)}(x)$, then we have the discrete correction equation:
\begin{subequations}
\label{eq:upwind_dg_dsa_correction}
\begin{align}
&\left(v_j (\BD_C- \frac{1}{2}\BD_J)+\BSigma_t\right)\left(\delta\bphi^{(l)}+3v_j\delta\BJ^{(l)}\right) = \BSigma_s \delta\bphi + \BSigma_s (\bphi^{(l,*)}-\bphi^{(l-1)}),\quad \text{if}\; v_j\geq 0,\\
&\left(v_j (\BD_C+\frac{1}{2}\BD_J)+\BSigma_t\right)\left(\delta\bphi^{(l)}+3v_j\delta\BJ^{(l)}\right) = \BSigma_s \delta\bphi + \BSigma_s (\bphi^{(l,*)}-\bphi^{(l-1)}),\quad \text{if}\; v_j<0.
\end{align}
\end{subequations}
Taking the discrete zero-th and first order moment of \eqref{eq:upwind_dg_dsa_correction} in the angular space, we obtain 
\begin{subequations}
\label{eq:dsa_first_moment}
\begin{align}
    &\BD_C\delta\BJ^{(l)} - (\sum_{v_j>0} \omega_jv_j) \BD_J\delta\brho^{(l)}+\BSigma_a \delta\bphi^{(l)} = \BSigma_s (\bphi^{(l,*)}-\brho^{(l-1)}),\label{eq:dsa_rho}\\
    (\sum_{j=1}^{N_v}\omega_j v_j^2)&\BD_C\delta\bphi^{(l)} +\left(\BSigma_t- 3(\sum_{v_j>0}\omega_j v_j^3)\BD_J\right)\delta\BJ^{(l)} = 0.\label{eq:dsa_J}
\end{align}
\end{subequations}
Equation \eqref{eq:dsa_J} yields
\begin{equation}
\delta\BJ^{(l)}=-(\sum_{j=1}^{N_v}\omega_j v_j^2)\left(\BSigma_t- 3(\sum_{v_j>0}\omega_j v_j^3)\BD_J\right)^{-1}\BD_C\delta\brho^{(l)}.
\label{eq:dsa_J_in_rho}
\end{equation} 
Substituting \eqref{eq:dsa_J_in_rho} into \eqref{eq:dsa_rho} to eliminate $\delta\BJ^{(l)}$, we obtain a fully consistent  DSA. If dropping the $\BD_J$ term in \eqref{eq:dsa_J_in_rho} before substituting it into \eqref{eq:dsa_rho}, we obtain the partially consistent DSA in \cite{wareing1993new}.

%%%%%%%%%%%%%%%%%%%%%%%%%%%%%%%%%%%%%%%%%%%%%%%%%%%%%%%%%%%%
\section{Modified Gram-Schmidt procedure with a truncation step\label{sec:mgs-procedure}}

\begin{algorithm}[H]
\caption{\pzc{Modified Gram-Schmidt procedure with a truncation step. \label{alg:mght}} }
\begin{algorithmic}[1]
\STATE{\pzc{Given a set of orthornormal basis $\BU_{r,1},\dots,\BU_{r,r}\in\mathbb{R}^{N}$ forming the columns of the matrix $\BU_{r}=(\BU_{r,1},\dots,\BU_{r,r})\in\mathbb{R}^{N\times r}$, a vector $\bpsi\in\mathbb{R}^N$ and the truncation tolerance $\epsilon_{\textrm{QR}}=10^{-13}$.}}
\STATE{\pzc{Define $\bu=\bpsi.$}}
\FOR{\pzc{$i=1:r$}}
    \STATE{\pzc{$R=\BU_{r,i}^T\bu$.}}
    \STATE{\pzc{$\bu:=\bu-R\BU_{r,i}$.}}
\ENDFOR
\IF{$||\bu||>\epsilon_{\textrm{QR}}$}
    \STATE{\pzc{Normalize the new basis as $\bu:=\bu/||\bu||$.}}
    \STATE{\pzc{Update the orthornormal basis as $\BU_{r+1}=(\BU_r,\bu)$.}}
\ELSE
    \STATE{\pzc{Do not update $\BU_r$.}}
\ENDIF
\end{algorithmic}
\end{algorithm}

%%%%%%%%%%%%%%%%%%%%%%%%%%%%%%%%%%%%%%%%%%%%%%%%%%%%%%%%%%%%
\bibliographystyle{elsarticle-num} 
\bibliography{ref}
\end{document}